\documentclass{amsart}

\usepackage{times, enumerate}
\usepackage[T1]{fontenc}
\usepackage[all,pdf]{xy}

\usepackage[colorlinks=true]{hyperref}
\hypersetup{urlcolor=blue, citecolor=red}

\newcommand\End{\operatorname{End}}

\usepackage[normalem]{ulem}
\renewcommand\sout{\bgroup\markoverwith
{\textcolor{red}{\rule[0.7ex]{3pt}{1.4pt}}}\ULon}

\usepackage{times, amsmath, amsbsy, amssymb, verbatim,amscd}
\usepackage{color,graphicx,fancybox} 
\DeclareGraphicsExtensions{.jpg,.pdf,.png,.ps,.eps}

\RequirePackage{color}
\definecolor{DarkBlue}{rgb}{0.0,0.0,0.55}
\definecolor{LightYellow}{rgb}{1.0,1.0,0.8}
\definecolor{mygreen}{rgb}{0.0,0.501,0.0}
\definecolor{mybrown}{rgb}{0.647,0.164,0.164}
\definecolor{lightsteelblue3}{rgb}{0.634,0.710,0.804}
\definecolor{gray}{rgb}{0.2,0.2,0.2}

\newcommand\CC{\mathbb C}
\newcommand\KK{\mathbb K}

\newcommand\NN{\mathbb N}

\newcommand\RR{\mathbb R}

\newcommand\ZZ{\mathbb Z}

\newcommand{\maA}{\mathcal A}
\newcommand{\maB}{\mathcal B}
\newcommand{\maC}{\mathcal C}

\newcommand{\maG}{\mathcal G}
\newcommand{\maH}{\mathcal H}

\newcommand{\maL}{\mathcal L}

\newcommand{\maP}{\mathcal P}

\newcommand{\maS}{\mathcal S}

\newcommand{\maU}{\mathcal U}
\newcommand{\maV}{\mathcal V}

\newcommand{\pa}{\partial}

\newcommand\ede{ \, := \, }
\newcommand\seq{ \, = \, }

\newcommand{\CI}{\mathcal C^{\infty}}
\newcommand\CIc{\mathcal{C}^\infty_c}

\newcommand\GR{{(\maG, d, r, i, \mu)}}
\newcommand\ver{\operatorname{\mathrm{ver}}}
\newcommand\units{{\maG^{(0)}}}
\newcommand\Cio{{\maC^{\infty,0}}}
\newcommand\oR{{\overline{\RR}}}

\newcommand\Diff{\operatorname{Diff}}

\newcommand\maSS{\maS_{\tilde \gamma, \tilde \gamma-1,
\tilde \gamma' + 2, \tilde \gamma' + 1}}

\newtheorem{theorem}{Theorem}[section]
\newtheorem{proposition}[theorem]{Proposition}
\newtheorem{corollary}[theorem]{Corollary}

\newtheorem{lemma}[theorem]{Lemma}
\newtheorem{problem}[theorem]{Problem}

\theoremstyle{definition}

\newtheorem{definition}[theorem]{Definition}

\newtheorem{remark}[theorem]{Remark}
\newtheorem{example}[theorem]{Example}
\newtheorem{examples}[theorem]{Examples}

\title[Schr\"odinger operators]{Schr\"odinger operators
with non-integer power-law potentials and Lie-Rinehart algebras}

\author[I. Beschastnyi]{Ivan Beschastnyi} \address{Université Côte d'Azur\\
  Centre Inria d'Université Côte d'Azur, 2004, route des Lucioles, 06902,
  Sophia Antipolis, France}
\email{ivan.beschastnyi@inria.fr}

\author[C. Carvalho]{Catarina Carvalho} \address{{Centro de Análise Matemática,
Geometria e Sistemas Dinâmicos, Dep. Matemática, Instituto Superior Técnico,
University of Lisbon, Av. Rovisco Pais 1, 1049-001 Lisbon, Portugal
}}
\email{catarinaccarvalho@tecnico.ulisboa.pt}

\author[V. Nistor]{Victor Nistor} \address{
Universit\'{e} de Lorraine, CNRS, IECL, F-57000 Metz, France}
\email{victor.nistor@univ-lorraine.fr}

\author[Y. Qiao]{Yu Qiao} \address{School of Mathematics and
  Statistics, Shaanxi Normal University, Xi'an, 710119,
China} \email{yqiao@snnu.edu.cn}

\thanks{I.B. was supported by the French government through the France
2030 investment plan managed by the National Research Agency (ANR), as part of
the Initiative of Excellence Université Côte d'Azur under reference number
ANR-15-IDEX-01, C.C. was supported by FCT/Portugal through project UIDB/04459/2020
with DOI identifier 10-54499/UIDP/04459/2020. V.N. was supported by ANR grant OpART
ANR-23-CE40-0016, Y.Q. was supported by NSFC (11971282, 12271323).}

\begin{document}

\maketitle

\begin{abstract}
    We study Schr\"odinger operators $H:= -\Delta + V$ with potentials $V$ that
    have power-law growth (not necessarily polynomial) 
    at 0 and at $\infty$ using methods of Lie
    theory (Lie-Rinehart algebras) and microlocal analysis. More precisely, 
    we show that $H$ is ``generated'' in a certain sense by an explicit 
    Lie-Rinehart algebra. This allows then to construct a suitable (microlocal) 
    calculus of pseudodifferential operators that provides further properties 
    of $H$. Classically, this microlocal analysis method was used to study $H$ 
    when the power-laws describing the potential $V$ have integer exponents.
    Thus, the main point of this paper is that this integrality condition
    on the exponents is not really necessary for the microlocal analysis 
    method to work. While we consider potentials 
    following (possibly non-integer) power-laws both at the 
    origin and at infinity, our results extend right away to potentials having
    power-law singularities at several points. The extension of the
    classical microlocal analysis results to potentials with non-integer 
    power-laws is achieved by considering the setting of Lie-Rinehart algebras 
    and of the continuous family groupoids integrating them. (The classical
    case relies instead on Lie algebroids and Lie groupoids.)
\end{abstract}

\tableofcontents

\section{Introduction}
\label{sec:intro}

We study the Schr\"odinger operator with \emph{power-law potential:}
\begin{equation}\label{eq.schrV}
    \begin{cases}
    \, H \ede -\Delta + V\,, & \\
    \ V(x) \sim |x|^{-a}\,, & \mbox{as } x \to 0\,, \ \mbox{ and}\\
    \ V(x) \sim |x|^{a'}\,, & \mbox{as } x \to \infty\,, \ a, a' \in \RR\,,
    \end{cases}
\end{equation}
using Lie-Rinehart algebras and microlocal analysis.
(See Section \ref{ssec:schroedinger} for precise definitions).
This type of potentials arises often in Physics, where $a$ is referred to 
as the ``hardness'' of an effective potential. In classical physics, this
type of potentials models non-hard collisions of particles in real systems, such
as colloids~\cite{Inverse_power_physics}. There exist many classical mathematical 
works dealing with Schr\"odinger operators with power-law potentials \emph{if $a$ 
and $a'$ are integers.} The main goal of
this paper is to show how the framework of \emph{Lie-Rinehart algebras}
and of \emph{continuous family groupoids} can be used to extend some of
these classical methods (i.e. for $a$ and $a'$ integers) to the general
case $a, a' \in \RR$ not necessarily integers.

Here is a concrete description of our results for the Schr\"{o}dinger
operator $H := -\Delta + V$ defined in Equation \eqref{eq.schrV}. Let $M = S^{n-1}
\times [0, \infty]$ (where $S^{k} \subset \RR^{k+1}$ is the
unit sphere) and let $\phi : M \to [0, \infty]$ be a smooth function such that:
$\phi(x) = |x|^{a}$ if $|x|$ is small, $\phi(x) = |x|^{-a'}$
if $|x|$ is large, and $\phi > 0$ on the rest of $M$.
\begin{enumerate}[(i)]
    \item A first result is to construct an algebra $\maA_{\maS}$ of
    functions
    \begin{equation}
        \CIc(M) \subset \maA_{\maS} \subset \maC(M)
    \end{equation}
    and a Lie algebra of vector fields $\maL_{\maS} \subset \CI(M; TM)$
    acting by derivations on $\maA_{\maS}$ and such that
    \begin{equation}\label{eq.belonging1}
        \phi (-\Delta + V) \in \Diff(\maA_{\maS}, \maL_{\maS})
    \end{equation}
    (that is, $\phi (-\Delta + V)$ belongs to the algebra of differential
    operators on $M$ generated by $\maA_{\maS}$ and
    $\maL_{\maS}$) see Proposition \ref{prop.inPsiS}.
    This type of results are known to be relevant for the study of the essential 
    spectrums of $\phi (-\Delta + V)$ and $-\Delta + V$ via general results for 
    operators in $\Diff(\maA_{\maS}, \maL_{\maS})$.

    \item A second result is to construct a suitable pseudodifferential
    calculus $\Psi(\maS)$ that ``quantizes'' the algebra
    $\Diff(\maA_{\maS}, \maL_{\maS})$ and which \emph{contains its parametrices.}
    Some canonical completions of $\Psi(\maS)$ will also contain the resolvents
    of the differential operators in $\Diff(\maA_{\maS}, \maL_{\maS})$.
    While $-\Delta + V$ is not in $\Psi(\maS)$ (it is $\phi (-\Delta + V)$
    who belongs to $\Psi(\maS)$), with some additional work, we can describe
    also the parametrices and resolvents of $-\Delta + V$ (Theorem \ref{thm.inPsiS}).

    \item Classically, the construction of pseudodifferential calculi is often
    achieved by ``integrating'' $(\maA_{\maS}, \maL_{\maS})$ to a Lie groupoid. 
    This is not possible in our case because the functions in $\maA_{\maS}$ and 
    the vector fields in $\maL_{\maS}$ are not smooth enough to come from a Lie algebroid
    (Lie algebroids are the infinitesimal forms of Lie groupoids).
    Nevertheless, we prove that the pair $(\maA_{\maS}, \maL_{\maS})$
    is a Lie-Rinehart algebra (Theorem \ref{thm.descr.A}), (a \emph{Lie-Rinehar}
    algebra is an algebraic generalization of a Lie algebroid).

    \item The natural object integrating a Lie-Rinehart algebra (when this is
    possible) is a CF-groupoid (the short hand for ``continuous
    family groupoid'' \cite{Paterson2000, Paterson2004}). The general
    problem of integrating a Lie-Rinehart algebra (formulated here) has not yet
    been solved. It will be
    studied in more detail in our forthcoming paper \cite{BCNQ_LR}.
    In this paper, we content to integrate the Lie-Rinehart algebra
    $(\maA_{\maS}, \maL_{\maS})$ to a CF-groupoid $\maS = \maSS$
    for some explicit parameters $\tilde \gamma$ and $\tilde \gamma'$ given in terms of 
    $a$ and $a'$ (and hence, in terms of the potentials). Our pseudodifferential 
    calculus is the pseudodifferential
    calculus associated to the CF-groupoid $\maS$ \cite{CNQ, LMN2000, Paterson2000}.
    As mentioned above, it can be shown to contain the parametrices of
    differential operators in $\Diff(\maA_{\maS}, \maL_{\maS})$.
\end{enumerate}

Technically, we first deal with the points (iii) and (iv) above
and then we use them to fully establish (i) and (ii). Of course, we establish 
many other intermediate results, some of which are of independent interest. 
For instance, a lot of the
technical work goes into the construction of the CF-groupoid $\maS = \maSS$
integrating our Lie-Rinehart algebra $(\maA_{\maS}, \maL_{\maS})$.
The extensive background material that we include is needed for this
construction. Nevertheless, we try to go as quickly as possible to the 
applications related to Schr\"{o}dinger operators, choosing to include 
further theoretical developments in the next paper in this series 
\cite{BCNQ_LR}.

In addition to allowing us to treat non-integer power-law singularities,
our approach (based on Lie-Rinehart algebras and CF-groupoids), allows us
to treat differential operators whose coefficients have lower regularity
near the boundary. This is critical in dealing with the Schr\"{o}dinger
operator \eqref{eq.schrV}.

Another example similar to the power-law Schr\"odinger operators of Equation
\eqref{eq.schrV}
involves Riemannian manifolds with boundary and metrics that, in a tubular
neighborhood of the boundary, are of the form
\begin{equation*}
    f \seq \frac{dx^2+g_x}{x^a}\,,
\end{equation*}
where $x$ is a boundary-defining function on our manifold
with boundary (i.e. $x \ge 0$, the boundary is given by $x = 0$,
and $dx$ does not vanish on the boundary) and $g_x$ is a family of
Riemannian metrics on the level sets of the defining function $x$.
As it is known, the case $a = 2$ corresponds to asymptotically hyperbolic
spaces. However, the Laplace-Beltrami operator on the same spaces with
fractional values of the parameter $a$ has been proposed recently for
studying acoustic modes of gas giants~\cite{De_verdiere+trelat=gas_giants}.
After a change of variables, these metrics can be put in the form
\begin{equation*}
    f \seq dx^2 +x^{-\beta} g_x.
\end{equation*}
Manifolds with metrics of this type are called \emph{Grushin manifolds.}
They have been studied in the cases of non-integer real parameter $\beta$
in a number of works~\cite{ivan_closure, ivan_selfadjoint_alpha_grushin,
ugo+dario_grushin, ame_alpha_cylinder, eugenio_alpha_plane}.

\subsection*{Contents of the paper}
In Section \ref{sec.cfgrpds}, we review the notions of Lie-Rinehart algebras and
CF-groupoids. We also briefly discuss locally compact and Lie groupoids.
In particular, we introduce the Lie-Rinehart algebra, the differential operators,
and the pseudodifferential calculus associated to a CF-groupoid \cite{LMN2000,
Paterson2000, Paterson2004}. We also recall the concept of the universal enveloping
algebra of a Lie-Rinehart algebra and discuss its relation to the algebra of
differential operators associated to a CF-groupoid.
At the end of this section, we examine
one-dimensional examples in detail.
In Section \ref{sec:general_integration+Schroedinger}
we construct the CF-groupoids used to study Schr\"odinger operators.
Our construction relies on the one dimensional examples
of the previous section. It is an instance of ``integrating''
a suitable Lie-Rinehart algebra.
In the last section, we apply the results of the previous sections to study
Schr\"{o}dinger operators with (possibly non-integer) power-law singularities at 0
and at $\infty$.

\subsection*{Acknowledgements}
We thank Camille Laurent-Gengoux and Robert Yuncken for useful discussions.
Some results of this paper were announced in \cite{BCNQproc}.

\section{Continuous family groupoids and the associated Lie-Rinehart algebras}
\label{sec.cfgrpds}

In this section we review continuous family groupoids and their associated
Lie-Rinehart algebras \cite{LMN2000, Paterson2000, Paterson2004}. We also briefly
discuss locally compact and Lie groupoids.

\subsection{Locally compact groupoids}
\label{ssec.lcg}
Let us begin by recalling the definition of groupoids in general
\cite{CNQ, ConnesNCG, LMN2000, MoerdijkFolBook, Paterson2000, RenaultB}.
A {\em groupoid} is a small category with all morphisms invertible.
More precisely, we have the following concrete description of groupoids.

\begin{remark}\label{rem.eqiv.groupoid}
    A groupoid is determined by the following data:
    \begin{enumerate}
        \item a set $\maG$ of \emph{morphisms};
        \item a set of \emph{composable pairs of morphisms}, denoted by
        $\maG^{(2)} \subset \maG \times \maG$,
        together with a \emph{multiplication map}
        $\mu:\maG^{(2)} \rightarrow \maG$, also written $\mu(g, h)=gh$; and
        \item an \emph{inverse map} $\iota : \maG \rightarrow \maG$, also written
        $\iota(g) = g^{-1}$.
    \end{enumerate}
    This data is required to satisfy the following conditions:
    \begin{enumerate}
        \item if $(g_1, g_2) \in \maG^{(2)}$ and $(g_2,g_2)\in \maG^{(2)}$,
        then \begin{equation*}(g_1g_2)g_3=g_1(g_2g_3),\end{equation*}
        (in particular, $(g_1g_2,g_3),\, (g_1, g_2g_3) \in \maG^{(2)}$, i.e.
        these pairs are composable);
        \item for all $g \in \maG$, we have $(g,g^{-1}) \in \maG^{(2)}$ and
        if $(g,h)\in \maG^{(2)}$, then
        \begin{equation*}
            g^{-1}(gh) \seq h \ \mbox{ and }\  (gh)h^{-1} \seq g\,.
        \end{equation*}
    \end{enumerate}
    We then define the {\em domain} and {\em range} maps
    $d: \maG \rightarrow \units$, $r: \maG \rightarrow \units$ by
    \begin{equation*}
        d(g) \ede g^{-1}g  \ \mbox{ and }\   r(g) \ede gg^{-1},
    \end{equation*}
    where the {\em unit space} $\units$ is defined to be $d(\maG)=r(\maG)$.
\end{remark}

Usually, we will write $\GR$ or, simply, $\maG$, to denote a groupoid.

A {\em locally compact groupoid} is a groupoid $\GR$ that is second countable,
locally compact space\footnote{Our locally compact spaces are not required
to be Hausdorff; some authors use the terminology ``locally quasi-compact''
spaces for these spaces.}
with the space of units $\units$ a closed, Hausdorff subset such that all
structural maps are continuous.
In the sequel, for subsets $A, B \subset \units$, we denote by $\maG_A:=d^{-1}(A)$,
$\maG^B:=r^{-1}(B)$, and $\maG_A^B:=d^{-1}(A)\cap r^{-1}(B)$.
In particular, if $x\in \units$, then $\maG_x$ and $\maG^x$ are the
\emph{$d$-fiber} and \emph{$r$-fiber at $x$,} respectively.
The set $\maG_A^A$ is also a groupoid, called the {\em reduction} of $\maG$ to $A$.
Moreover, if $\maG_A^A=\maG_A=\maG^A$, then $A$ is called {\em invariant}.
In this case, $\maG_A$ is also a groupoid.
For any $x\in \units$, $\maG_x^x$ is a group, called the {\em isotropy group} at $x$.

\subsection{Continuous family groupoids}
\label{ssec:cfg}
Our motivation for studying groupoids is that Lie groupoids have been
used to study certain partial differential equations with
\emph{singularities}, see
\cite{aln1, YunckenNAll, andr22, ASkandalis1, ivan_closure, BCNQproc,
CNQ, ConnesNCG, kottkeInvent, Mazzeo91, Savin05, vanErpYuncken}
and the references therein. As explained in the introduction, certain differential
operators (such as our Schr\"{o}dinger operators with non-integer
power-law singularities) cannot be studied using Lie groupoids. We prove, however,
that they can be studied using continuous family groupoids. We will do that
for Schr\"odinger operators in the last subsection of this paper.

Continuous family groupoids (\emph{CF-groupoids,} for
short) were introduced by Paterson
\cite{Paterson2000},
who has studied index problems in this setting \cite{Paterson2004}.
Broadly speaking, a continuous family groupoid
is such that all $d$-fibers $\maG_x:= d^{-1}(x)$ and all $r$-fibers
$\maG^{x}:=r^{-1}(x)$ are smooth manifolds for $x\in \units$, and these fibers
vary continuously on $x$, forming a \emph{continuous family of manifolds} over
the units.

We next briefly review the necessary material on continuous families
of manifolds and on continuous family groupoids from \cite{Paterson2000}.

\subsubsection{Fiber spaces}
Let $X$ and $Y$ be topological spaces and $p: X \rightarrow Y$ be
a continuous, open surjection. The pair $(X,p)$ is then called a fiber space
over $Y$ with fibers $X^y:=p^{-1}(y)$ for $y \in Y$. Let $(X_1, p_1)$ and $(X_2,p_2)$
be fiber spaces over $Y_1$ and $Y_2$, respectively, and
let $q: Y_1 \rightarrow Y_2$ be a continuous map. A continuous map
$f:X_1 \rightarrow X_2 $ will be called
\emph{fiber-preserving with respect to $q$} if $p_2\circ f= q \circ p_1$. In case
$Y_1 \subset Y_2$ and $q : Y_1 \to Y_2$ is the inclusion, a fiber-preserving map 
$f : X_1 \to X_2$
with respect to the inclusion $q$ will be called simply \emph{fiber-preserving.}

Let $(X, p)$ be a fiber space over $Y$ (so $p : X \to Y$).
Let $Z$ be a topological space and $\rho: Z \rightarrow Y$ be a continuous map.
The {\em pull-back} $(\rho^{-1}(X), \pi_2)$ of $(Y, p)$ over $Z$ (along the map
$\rho$) is defined to be
\begin{equation}\label{eq.def.pull-back}
    \rho^{-1} (X) \ede \{(x, z)\in X \times Z\mid \rho(z)=p(x) \in Y \}
    \subset X \times Z\,,
\end{equation}
with fibration map $\pi_{2}(x, z) \ede z$ and the topology induced from
the product $X \times Z$. It is a fiber space over $Z$. Let $(W, q)$ be a
fiber space over $Z$. Then a fiber-preserving map $f : W \to X$ over $\rho$ is
the same thing as a fiber preserving map $f_\rho : W \to \rho^{-1}(X)$ (as $Z$-fiber
spaces), where $f_\rho(w) = (f(w), q(w))$. Indeed, $f$ is fiber-preserving over
$\rho$ if, and only if, $\rho(q(w)) = p(f(w))$, which in turn is equivalent to
$f_\rho(w) \in \rho^{-1}(X)$. Moreover, $f$ is continuous if, and only if,
$f_\rho$ is continuous. This allows to reduce the study of ``fiber-preserving maps
over a map $\rho$'' to simply ``fiber-preserving
maps'' (i.e. ``fiber-preserving maps over the identity'').

\subsubsection{Continuous families of manifolds and $\Cio$ maps}
The main difference between Lie groupoids and CF-groupoids is that, in the later,
smooth manifolds are replaced by ``continuous families of manifolds,'' a concept
that we introduce in this section. To this end, we now look at a particular
case of fiber-preserving maps.

\begin{remark}
    \label{rem.cinf.fam}
    Let $Y$ be a topological space, $A$ and $B$ be open subsets of $ \RR^{k} \times Y$
    for some $k \geqslant 1$, and $\pi_2$ be the canonical projection from
    $\RR^{k} \times Y$ to $Y$. Then $(A, \pi_{2})$ and $(B, \pi_{2})$ are fiber spaces
    over subsets of $Y$ and a continuous map $f : A \rightarrow B$ is fiber-preserving if, 
    and only
    if, for every $y \in A$, we have $\pi_2(f(y))= \pi_2(y)$. The definition implies that
    $\pi_2(A) \subset \pi_2(B)$. Let us consider a \emph{fiber-preserving function}
    $f: A \rightarrow B$ and let $U_i \subset \RR^k$ and $V_i \subset Y$, $i=1,2$, be
    open subsets such that
    \begin{equation*}
        U_1 \times V_1 \subset A\,,\ U_2\times V_2 \subset B\,,\ \mbox{ and }\
        f(U_1\times V_1) \subset U_2\times V_2\,.
    \end{equation*}
    Then, for each $y \in V_1$, there exists a unique function $f_y : U_1 \to U_2$
    such that
    \begin{equation}\label{eq.def.fiber.fy}
        f(x, y) \seq (f_y(x), y)\in U_2 \times V_2\,, \quad \mbox{where }\
        x \in U_1 \mbox{ and } y \in V_1\,.
    \end{equation}
    We let $\CI(U_1, U_2)$ denote the set of \emph{smooth} functions $U_1 \to U_2$
    endowed with the topology of uniform convergence on compacts for all partial
    derivatives.
\end{remark}

We shall often use the notation and concepts introduced in the above
remark, without further discussion.

\begin{definition}[Connes \cite{ConnesFol}, Paterson \cite{Paterson2000}]
    \label{def.long.smooth.f1}
    Let $A, B \subset \RR^{k}\times Y$ with the induced fiber
    spaces structures over (subsets of) $Y$, as in Remark \ref{rem.cinf.fam},
    whose notation we continue to use. A fiber-preserving function
    $f: A \rightarrow B$ is called a {\em $\Cio$-function} (or \emph{longitudinally
    smooth}) if, whenever $U_i \subset \RR^k$ and $V_i \subset Y$, $i=1,2$, are
    open subsets such that $U_1\times V_1 \subset A$, \ $U_2\times V_2 \subset B,$
    and $f(U_1\times V_1) \subset U_2\times V_2$, the induced map $y \rightarrow f_y$
    is continuous from $V_1$ to $\maC^{\infty}(U_1, U_2)$
    (see Equation \ref{eq.def.fiber.fy} for the definition of $f_{y}$).
    We let $\maC^{\infty, 0}(A,B)$ denote the set of $\Cio$-functions from $A$ to $B$.
    A $\Cio$ fiber-preserving homeomorphism $f : A \to B$ such that $f^{-1}$ exists and
    $f^{-1} \in \maC^{\infty, 0}(B,A)$ will be called a \emph{$\Cio$-homeomorphism.}
\end{definition}

This definition allows us now to introduce continuous families
of manifolds, as in \cite{Paterson2000}.

\begin{definition}\label{def.cont.fam.man}
    A fiber space $(X,p)$ over $Y$ is called {\em a continuous family
    of smooth manifolds}
    over $Y$ if it is paracompact and there exists a set of coordinate
    charts
    $\{(U_\tau, \varphi_\tau) \}_{\tau\in I}$, where
    every $U_\tau$ is an open subset of $X$ and
    $\cup_{\tau \in I} U_{\tau} = X$,
    such that
    \begin{enumerate}
        \item for every index $\tau$, the map $\varphi_\tau$ is a fiber-preserving
        homeomorphism from $U_\tau$ to an open subset of $ \RR^k \times Y$
        (recall that fiber preserving means
        $\pi_2\circ \varphi_{\tau} = p \vert _{U_\tau}$);

        \item for every pair of indices $\tau$ and $\beta$, the composition map
        \begin{equation*}\varphi_\beta \circ \varphi_\tau^{-1} :
            \varphi_\tau(U_\tau \cap U_\beta) \to
        \varphi_\beta(U_\tau \cap U_\beta)\end{equation*}
        is a $\Cio$-homeomorphism (see Definition \ref{def.long.smooth.f1}).

        \item The fiber $p^{-1}(y) \subset X$ of $p$ is Hausdorff
        for every $y \in Y$.
    \end{enumerate}
\end{definition}

It follows right away from this definition that, if $(X, p)$ is a continuous
family of smooth manifolds over $Y$, then, for every $y \in Y$, the space
$p^{-1}(y)$
has a natural smooth manifold structure in the usual sense. More generally,
if $q : Z \to Y$ is a continuous map and
$(X, p)$ is a continuous family of smooth manifolds over $Y$, then the
pull-back
$(q^{-1}(X), \pi_{2})$ is a continuous family of smooth manifolds over $Z$
\cite{Paterson2000} (see Equation \eqref{eq.def.pull-back} for the
definition of
$(q^{-1}(X), \pi_{2})$).

Recall also that our smooth manifolds \emph{do not have boundaries or corners.}
Thus a ``manifold with corners'' is not a ``smooth manifold''
in our sense. We shall need the following lemma.

\begin{lemma}\label{lemma.projection}
    Let $(X, p)$ be a continuous family of smooth manifolds over $Y$
    and let $f : X \to \CC$ be a continuous function such that
    $f(x_1) = f(x_2)$
    for all $x_1, x_2 \in X$ with $p(x_1) = p(x_2)$. Then there exists a
    continuous function $f_0 : Y \to \CC$ such that $f = f_{0} \circ p$.
    The converse is also true (and trivial). Similarly, let $E \to Y$
    be a (locally trivial) vector bundle and $\xi$ be a continuous section of
    $p^{*}(E)$. There exists a continuous section
    $\xi_{0} : Y \to E$ such that $\xi = \xi_{0} \circ p$ if,
    and only if, $\xi(x_1) = \xi(x_2)$
    for all $x_1, x_2 \in X$ with $p(x_1) = p(x_2)$.
\end{lemma}

\begin{proof}
    The existence of a function $f_0$ with these properties is immediate. We
    only need to prove that it is continuous. Indeed, by the definition of a
    continuous family of manifolds, it follows that $f_0$ is locally continuous.
    But a locally continuous function is continuous. The converse is
    trivial since $p$ is continuous and hence $f := f_0 \circ p$ is also
    continuous
    and has the desired properties. To prove the extension of the
    statement to vector bundles, we first notice that it is a
    direct consequence of the first part if $E = Y \times \RR^{k}$
    (a trivial vector bundle). By localizing, we can then reduce the
    general case to that of a trivial vector bundle.
\end{proof}

We now extend the definition of $\Cio$-maps between subsets of
$\RR^{k}\times Y$,
Definition \ref{def.long.smooth.f1} to arbitrary fiber spaces.

\begin{definition}\label{def.smooth.funct}
To define a \emph{$\Cio$-function} between two continuous families of
manifolds $(X_1, p_1)$ and $(X_2,p_2)$, it is convenient to distinguish
two cases.
\begin{enumerate}
    \item If two continuous families $(X_1, p_1)$ and $(X_2,p_2)$
    of manifolds are over the same base $Y$, a $\Cio$-function from
    $(X_1, p_1)$ to $(X_2,p_2)$ \emph{(covering the identity of $Y$)}
    is simply a function $f : X_{1} \to X_{2}$ such that,
    whenever $(U,\varphi)$ and $(U', \varphi')$
    are local charts for $X_1$ and $X_2$ respectively with $p_1(U) \subset
    p_2(U')$ and $f(U)\subset U'$ (see Definition \ref{def.cont.fam.man}),
    then $\varphi' \circ f \circ \varphi^{-1} \in \maC^{\infty, 0}
    (\varphi(U), \varphi'(U'))$. If this is the case, we write
    $f \in \maC^{\infty, 0}(X_1, X_2)$.

    \item Let $(X_1, p_1)$ and $(X_2,p_2)$ be two continuous families of manifolds over
    $Y_1$ and $Y_2$ respectively. Let $q: Y_1 \rightarrow Y_2$ be a continuous map
    and $f:X_1 \rightarrow X_2 $ be a continuous fiber-preserving map with respect to $q$,
    that is, $p_2\circ f= q \circ p_1$. Recall the  pull-back $q^{-1}(X_2)$, Equation
    \eqref{eq.def.pull-back}. Then $f$ defines a continuous, fiber preserving map
    $f_q : X_1 \to q^{-1}(X_2) \subset X_2 \times Y_1$ by the formula
    $f_q(x) = (f(x), p_1(x))$.
    Notice that both $X_1$ and $q^{-1}(X_2)$ are
    families of smooth manifolds over the same base $Y_1$. We say that $f$ is a
    $\maC^{\infty, 0}$-function \emph{(covering $q : Y_1 \to Y_2$)} if $f_q$ is a
    $\maC^{\infty, 0}$-function \emph{(covering the identity map of $Y_1$)}.
\end{enumerate}
\end{definition}

It is not hard to show that the composite of two $\Cio$-maps of continuous
families of manifolds is still a $\Cio$-map. A $\Cio$-map is, in particular,
fiber preserving and smooth on each fiber.

\begin{remark}
\label{rem.vert.tb}
As in the case of usual (single) smooth manifolds, to a continuous family of
manifolds $(X, p)$ over $Y$, we can associate in a functorial way its (vertical)
tangent and cotangent bundles $T_{\ver}X$ and $T_{\ver}^{*}X$, see
\cite{Paterson2000}.
\end{remark}

\subsubsection{Fibered products and continuous family groupoids}
Suppose that $(X_1, p_1)$ and $(X_2, p_2)$ are two continuous families of manifolds over $Y$.
The {\em fibered product} of $(X_1, p_1)$ and $(X_2, p_2)$ is defined to be
\begin{equation}\label{def.fib.prod}
    (X_1, p_1) * (X_2, p_2) \ede \{(x_1, x_2)\in X_1 \times X_2 \,
    \vert \, p_1(x_1)=p_2(x_2)\} \subset X_1 \times X_2\,,
\end{equation}
with fiber map $p : (X_1, p_1) * (X_2, p_2) \to Y$ given by
$p(x_1,x_2)=p_1(x_1)=p_2(x_2)$. We shall usually write $X_1 * X_2$ instead of
$(X_1, p_1) * (X_2, p_2)$ when the maps $p_1$ and $p_2$ are obvious.
The fiber space $(X_1 *  X_2, p)$ is a continuous
family of smooth manifolds over $Y$ with respect to the subspace topology. Moreover,
the fibered product $X_1 *  X_2 $ can naturally be regarded as a continuous family of manifolds
over $X_1$ and $X_2$ respectively. For instance, the pullback of $(X_2,p_2)$ along the map
$p_1: X_1 \rightarrow Y$ gives the continuous family $(p_1^{-1}(X_2), \pi_1) = (X_1 * X_2, \pi_1)$
of manifolds over $X_1$, which is illustrated in the following commutative
diagram:
\begin{equation*}
\xymatrix{
    X_1*X_2=p_1^{-1}(X_2) & X_2\\
    X_1 & \ Y
    \ar^{\qquad \quad\pi_2}"1,1";"1,2"
    \ar^{\pi_1}"1,1";"2,1"
    \ar^{p_2}"1,2";"2,2"
    \ar^{p_1}"2,1";"2,2"
}
\end{equation*}

Here we switch the first and second components in the definition of the pullback 
(see Equation (\ref{eq.def.pull-back})) to match Equation (\ref{def.fib.prod}).
Likewise, we are able to obtain a continuous family $(X_1* X_2, \pi_2)$ of manifolds
over $X_2$. We are now in position to recall the definition of continuous family
groupoids from \cite{Paterson2000}.

\begin{definition}\label{def.CF.groupoid}
    A groupoid $\GR$ is called a \emph{continuous family groupoid}
    (\textit{CF-groupoid }for short) if the following conditions hold:
    \begin{enumerate}
        \item The set $\units$ is locally compact Hausdorff and both
        $(\maG, d)$ and $(\maG,r)$ are continuous families of manifolds over $\units$;
        \item the inverse map $\iota :
        (\maG,d) \rightarrow (\maG,r)$, where $i(g)=g^{-1}$, is a
        $\maC^{\infty, 0}$-homeomorphism of continuous families of manifolds;
        \item the multiplication map $\mu: \big((\maG,d)*(\maG,r), \pi_1\big)
        \rightarrow (\maG,r)$ is a $\Cio$-map of continuous
        families of manifolds with respect to the range map $r$:
        \begin{equation*}
            \xymatrix{
            \maG * \maG & {\ \maG \ } \\
            {\ \maG \ } & {\ \units }
            \ar^{\mu}"1,1";"1,2"
            \ar^{\pi_1}"1,1";"2,1"
            \ar^{r}"1,2";"2,2"
            \ar^{r}"2,1";"2,2"
            }
        \end{equation*}
    \end{enumerate}
    If $\maG$ and $\maG'$ are two CF-groupoids, a \emph{CF-groupoids isomorhism}
    $\Phi : \maG \to \maG'$ is an isomorphism of groupoids that is also
    a $\Cio$-diffeomorphism.
\end{definition}

Roughly speaking, a continuous family groupoid
$\GR$ is such that all $d$-fibers $\maG_x:= d^{-1}(x)$ (respectively $r$-fibers
$\maG^{x}:=r^{-1}(x)$) are smooth manifolds for $x\in \units$, and these fibers
are required to vary continuously on $x$. Thus $\maG_{x}$ does not have
boundary or corners. The structural maps are required to be smooth along the fibers
and continuous along the units. An equivalent definition of continuous family groupoids
using coordinate charts is given in \cite{LMN2000}.

\begin{examples} \label{examples}
    The following examples will play an important role.
    \begin{enumerate}
        \item \label{item.pair}
        Let $M_{1}$ be a smooth manifold and
        $\maP := M_{1} \times M_{1}$ with $d$ the second projection and $r$ the
        first projection. Then
        $\maP$ is a Lie groupoid, and hence a CF-groupoid as well, called
        the \emph{pair groupoid}. However, if $M_{1}$ is a manifold with
        a non-empty boundary, then the pair groupoid $M_{1}^{2}$ is not a
        CF-groupoid.

        \item If, moreover, $Y$ is a topological space, then
        $Y \times \maP$ will be a CF-groupoid. It will be a Lie groupoid
        precisely if $Y$ is a manifold (possibly with boundary) with
        $d(y, m_{1}, m_{2}) := (y, m_{2})$ and
        $r(y, m_{1}, m_{2}) := (y, m_{1})$.

        \item
        If $\maG$ is a groupoid and $\maG$ and $\units$ are manifolds
        (possibly with corners), all the structural maps $d, r, \mu, \iota, u$
        are smooth, and $d$ is a submersion, then $\maG$ is a \emph{Lie groupoid.}
        (Recall that the condition that $d$ be a submersion implies that the
        fibers $\maG_{x} := d^{-1}(x)$ are smooth, boundaryless manifolds.)
        All Lie groupoids are, in particular, CF-groupoids.

        \item \label{item.cor}
        Let $G$ be a Lie group acting on a CF-groupoid $\maH$
        continuously via $\Cio$-isomorphisms. Then we form the \emph{action} (or
        \emph{semi-direct product}) groupoid $\maH \rtimes G$ with the
        same units as $\maH$ and with structural maps $d(h, g) = d(h)$,
        $r(h, g) = g(r(h))$, and product $(h, g)(h', g') = (g^{-1}(h')h, gg')$.
        Then $\maH \rtimes G$ is a CF-groupoid, \cite[Proposition 5]{Paterson2000}

        \item \label{item.thm} 
        Let $\maG$ be a continuous family groupoid with units $\units$
        and $(X,p)$ be a fiber space over $\units$. We consider the fiber space 
        $(\maG, d)$ and the fiber product $\maG* X$, see Equation 
        \eqref{def.fib.prod}. A \emph{continuous action} of $\maG$ on $X$ is a continuous map 
        $\rho: \maG* X \rightarrow X$ ($\rho(g,x)=gx$, for short) 
        that satisfies the following conditions:
        \begin{enumerate}
            \item $(g_1 g_2)y=g_1(g_2y)$, if $(g_{1}, g_{2}) \in \maG^{(2)}$
            (the set of composable
            arrows),
            \item $p(gy) = r(g)$, and
            \item $g^{-1}(gy)=y$,
        \end{enumerate}
        whenever they make sense. Given this data, we can form the \emph{semi-direct product}
        $X \rtimes \maG := (\maG,d)*(X,p)$. Then the semi-direct product $X \rtimes \maG := 
        \maG*X$ is a CF-groupoid with units $X$ \cite[Proposition 5]{Paterson2000}.
        If $\maG$ is a Lie groupoid, $X$ is a smooth 
        manifold, and the action is smooth, then the semi-direct product groupoid 
        $X \rtimes \maG := \maG * X$ is still a Lie groupoid \cite[Page 125]{MoerdijkFolBook}.
    \end{enumerate}
\end{examples}

We shall not use Lie groupoids in this paper for purposes other than to provide
examples, but the reader can find more information on Lie groupoids in
\cite{BCNQproc, ConnesNCG, MackenzieBook1, MackenzieBook1}.

\subsection{The Lie-Rinehart algebra associated to a CF-groupoid}
\label{ssec:lra_of_cfg}
Lie-Rinehart algebras constitute an algebraic generalisation of Lie algebroids.
Recall that a Lie algebroid  is given by a smooth vector bundle $A\to M$
with a Lie algebra structure on its space of smooth sections and a vector bundle
map $\rho:  A\to TM$ that satisfies a Leibniz rule. Before discussing Lie-Rinehart
algebras associated to continuous family groupoids, let us recall the definition
of Lie-Rinehart algebras~\cite{Rinehart} and give some examples.

\begin{definition}[\cite{Rinehart}]
    \label{def.LieRinehart}
    A (unital) \emph{Lie-Rinehart algebra} is a pair $(\maA, \maL)$ that consists
    of a unital commutative algebra $\maA$ and a Lie algebra $\maL$ with the
    following additional structures:
    \begin{enumerate}
    \item $\maL$ is a left $\maA$-module;
    \item $\maL$ acts by derivations on $\maA$ through
    an $\maA$-module Lie algebra homomorphism $\rho: \maL \rightarrow
    \operatorname{Der}(\maA)$
    \item The homomorphism $\rho$ satisfies the Leibnitz rule:
    \begin{equation*}[X, a Y] \seq \rho(X)(a) Y + a[X, Y]\end{equation*}
    for all $X, Y \in \maL, a, b \in \maA$.
    \end{enumerate}
    The map $\rho$ is called the anchor morphism.
    If $(\maA, \maL)$ and $(\maB, \mathcal{M})$ are two Lie-Rinehart
    algebras, a \emph{Lie-Rinehart isomorhism} $\Phi := (\Phi_{\maA}, \Phi_{\maL})
    : (\maA, \maL) \to (\maB, \mathcal{M})$ is a pair of an algebra and a Lie-algebra
    isomorphisms compatible with the respective actions.
\end{definition}

Lie-Rinehart algebras were considered also in 
\cite{Herz, HuebschmannPoisson, HuebschmannDuality, Huebschmann2004, CamilleBook, 
MoerdijkMrcunLieRine, Palais_Liering, Rinehart}. 
From now on everywhere in the text we assume that $\maA$ is an algebra
over $\mathbb{K}$, where $\mathbb{K}$ is either $\CC$ or $\RR$. The same
is true for all function spaces, which we assume to be either real or complex
valued. Generally, our vector bundles will be over $\RR$,
whereas all algebras (Lie or otherwise) will be over $\CC$.

Let $\GR$ be a CF-groupoid. We now introduce its associated Lie-Rinehart algebra
by analogy with the corresponding construction of the Lie algebroid of a Lie groupoid.
Let $\maG_x := d^{-1}(x)$, $x$ a unit of $\maG$ (that is, $x \in \units$). Then $\maG_x$
is endowed with the structure of \emph{a Hausdorff, smooth manifold,} by the
definition of a CF-groupoid (recall that our smooth manifolds do not have boundaries
or corners). Let $g \in \maG$, then, again by definition, the right multiplication by
$g$ induces a \emph{diffeomorphism}
\begin{equation}\label{eq.ref}
    \maG_{r(g)} \ni h \to hg \in \maG_{d(g)} \,,
\end{equation}
such that $R_g$ varies continuously with $g$ \cite{Paterson2000, Paterson2004}.
The same holds for left translations. 
Recall \cite{Paterson2000} that the tangent spaces to the fibers
$d^{-1}(x)$ of $d$ fit into a (continuous) vector bundle $T_{\ver}\maG$
determined by
\begin{equation} \label{eq.restr.cond}
    T \maG_x \seq T_{\ver}\maG\vert_{\maG_x}\,.
\end{equation}
The vector bundle $T_{\ver}\maG$ is called the \emph{$d$-vertical tangent bundle to $\maG$.}

\begin{definition}\label{def.LRalgebra}
    Let $\maG$ be a CF-groupoid. We let
    $\maA(\maG)$ consist of all the continuous functions $f: \maG \to \CC$
    with the following properties:
    \begin{enumerate}
        \item $f$ is $\Cio$ (recall that this means that it is continuous and its
        restriction to each $\maG_x := d^{-1}(x)$ is smooth).
        \item $f$ is right $\maG$-invariant (that is, for all $g,h \in \maG$,
        $d(h) = r(g)$, we have $f(hg) = f(h)$).
    \end{enumerate}
    Similarly, we let $\maL(\maG)$ consist of all the continuous sections $X$ of
    $T_{\ver}\maG$ with the following properties:
    \begin{enumerate}
        \item $X$ is $\Cio$ (or longitudinally smooth) in the following sense:
        it is a continuous section of $T_{\ver}\maG$ and its restriction to each
        $\maG_x := d^{-1}(x)$ is a smooth vector field.
        \item $X$ is right $\maG$-invariant (that is, for all $g,h \in \maG$,
        $d(h) = r(g)$, we have $X(hg) = g_*X(h)$).
    \end{enumerate}
    Then the pair $(\maA(\maG), \maL(\maG))$ is called
    \emph{the Lie-Rinehart algebra associated to $\maG$}.
\end{definition}

This definition is justified by the following lemma.

\begin{lemma} \label{lemma.LRA}
    The pair $(\maA(\maG), \maL(\maG))$ in Definition \ref{def.LRalgebra} 
    is a Lie-Rinehart algebra.
\end{lemma}

\begin{proof}
    Both the longitudinal smoothness and the right invariance properties
    are stable under multiplication or under the Lie bracket,
    and hence $\maA(\maG)$ is a commutative algebra
    for the pointwise multiplication and $\maL(\maG)$ is a Lie algebra
    for the Lie bracket of vector fields. Let $f \in \maA(\maG)$
    and $X \in \maL(\maG)$. Then $fX \in \maL(\maG)$ and $X(f) \in \maA(\maG)$
    by the definitions of $\maA(\maG)$, of $\maL(\maG)$, and of a continuous family
    of smooth manifolds (in this case, applied to $d : \maG \to \units$).
    Therefore $\maL(\maG)$ is also a $\maA(\maG)$-module and acts
    on $\maA(\maG)$, by derivations. The Leibnitz rule is satisfied since the Lie bracket of
    vector fields satisfies it.
\end{proof}

\begin{definition}\label{def.algebroid}
    Let $\maG$ be a CF-groupoid and $T_{\ver}\maG$ be the vertical
    tangent bundle to $\maG$ (tangent to
    the fibers of $d : \maG \to \units$). Then
    \begin{equation*}
        \mathbf{A}(\maG) \ede T_{\ver}\maG\vert_{\units}
    \end{equation*}
    is called the \emph{algebroid associated to} $\maG$.
    We shall also say that $\maG$ \emph{integrates} any Lie-Rinehart
    algebra isomorphic to $(\maA(\maG), \maL(\maG))$.
\end{definition}

Clearly, the algebroids of two isomorphic CF-groupoids are also isomorphic. The same is 
true for their associated Lie-Rinehart algebras.

\begin{remark}
    Let us assume that $\maG$ is a Lie groupoid with units $\units$.
    Recall that any Lie groupoid has an associated Lie algebroid $\mathbf{A}(\maG)$.
    In turn, to $\mathbf{A}(\maG)$, one can associate the Lie-Rinehart algebra
    $(\CI(\units), \CI(\units, \mathbf{A}(\maG)))$. This is \emph{not},
    however, the Lie-Rinehart algebra associated to $\maG$ when
    regarded as a CF-groupoid (although the algebroid
    $\mathbf{A}(\maG)$ will be the same). The reason is that, while 
    $\CI(\units)$ is always contained in $\maA(\maG)$, the two algebras 
    are usually different. The same comment applies to $\CI(\units, \mathbf{A}(\maG))$
    and $\maL(\maG)$: the second space always contains the first one, but 
    is usually larger. See Remark \ref{rem.inj.restr}.
    below.
\end{remark}

Recall that $\maC(K)$ denotes the space of continuous functions on some
locally compact space $K$ and that $\maC(K; E)$ denotes the space of continuous
sections of some vector bundle $E \to K$.

\begin{remark}\label{rem.inj.restr} We use the notation of
    Definition \ref{def.algebroid},
    then $\mathbf{A}(\maG)$ is a continuous vector bundle on $\units$
    with global (continuous) sections denoted $\maC(\units;
    \mathbf{A}(\maG))$. The restriction maps
    \begin{equation*}
        \begin{gathered}
            \maA(\maG) \ni  f \to f\vert_{\units} \in \maC(\units) \quad \mbox{and}\\
            \maL(\maG) \ni  X \to X\vert_{\units} \in \maC(\units;
            \mathbf{A}(\maG))
        \end{gathered}
    \end{equation*}
    are injective and we have a canonical isomorphism
    $T_{\ver}\maG \simeq r^*\mathbf{A}(\maG).$

\end{remark}

The algebroid $\mathbf{A}(\maG)$ associated to the CF-groupoid $\maG$ was introduced
by Paterson \cite{Paterson2000}. The algebroid $\mathbf{A}(\maG)$
was used in \cite{Paterson2000} to construct a continuous Haar system on $\maG$
using the $1$-density bundle of $\mathbf{A}(\maG)$. Let us record the following simple
lemma for further use.

\begin{lemma} \label{lemma.quotient}
    Let $f : \maG \to \CC$ be a continuous, $\maG$-right invariant
    function (that is, $f(hg) = f(h)$, whenever $g, h \in \maG$ are
    such that $d(h)=r(g)$). Then there exists a continuous function
    $f_0 : \units \to \CC$ such that, for all $g \in \maG$,
    $f(g) = f_0(r(g))$. The converse is trivially true: any function
    of the form $f_0 \circ r$ is $\maG$-right invariant.
    A similar statement holds for sections of $T_{\ver}\maG$: if $X$
    is a right invariant vertical vector field on $\maG$, then there
    exists a section $Y$ of $\mathbf{A}(\maG)$ such that
    $r_{*}(X) = Y$.
\end{lemma}

\begin{proof}
    This is a consequence of the fact that $r : \maG \to \units$ is
    a continuous family of smooth manifolds, of Lemma
    \ref{lemma.projection}, and (for the last part) of the definition
    relation $T_{\ver}\maG \simeq r^{*}(\mathbf{A}(\maG))$.
\end{proof}

\begin{remark}\label{rem.new.r}
    We let the spaces of restriction to $\units$ be given by
    \begin{equation*} %\label{eq.restr}
        \begin{gathered}
            \maA(\maG)_0:= \{ f_0 \in \maC(\units) \mid f_0\circ r \in
            \maA(\maG) \} \quad \mbox{and}\\
            \maL(\maG)_0:= \{ X\ \in \maC(\units; \mathbf{A}(\maG))
            \mid X\circ r \in    \maL(\maG) \}.
        \end{gathered}
    \end{equation*}
    The restriction maps define isomorphisms $\maA(\maG) \simeq
    \maA(\maG)_0 \subset \maC(\units)$ and $\maL(\maG) \simeq
    \maL(\maG)_0 \subset \maC(\units, \mathbf{A}(\maG))$. In particular,
    $( \maA(\maG)_0, \maL(\maG)_0)$ is a Lie-Rinehart algebra.
\end{remark}

\subsection{Universal enveloping algebras, differential,
and pseudodifferential operators}
\label{ssec:universal_envelop}
From now on, $(\maA, \maL)$ will denote a unital
Lie-Rinehart algebra (Definition \ref{def.LieRinehart}) over $\KK = \RR$
or $\KK = \CC$.
We shall need the concept of the \emph{universal enveloping algebra}
$\maU(\maA, \maL)$ of $(\maA, \maL)$ \cite{CamilleBook,
MoerdijkMrcunLieRine, Rinehart}, which we recall now. First, we endow the
$\KK$-vector space $\maA\oplus \maL$ with the Lie algebra structure with bracket
\begin{equation*}\label{eq.bracket.1}
    \left[ (a, X), (b,Y)\right] \ede \left( X(b)-Y(a), [X,Y] \right)\,.
\end{equation*}
Denote by $\maU(\maA\oplus \maL)$ the \emph{universal enveloping
algebra} of the Lie algebra $\maA\oplus \maL$, which is isomorphic to the quotient
of the tensor algebra
\begin{equation*}\label{eq.def.TA}
    T(\maA\oplus \maL) \ede \mathbb{K} \oplus \left [\, \oplus_{n \ge 1}
    (\maA\oplus \maL)^{\otimes n}\, \right ] \supset
    (\maA\oplus \maL)^{\otimes 1} \ede \maA\oplus \maL
\end{equation*}
by the ideal generated by the elements of the form
$u\otimes v- v\otimes u - [u,v] \in T(\maA\oplus \maL)$,
$u,v \in \maA \oplus \maL$. Let $i :  \maA\oplus \maL \to
\maU(\maA\oplus \maL)$ be the canonical inclusion map above, and
$\tilde{\maU}(A\oplus L)$
be the subalgebra of $\maU(\maA\oplus \maL)$ generated by the image of $i$.
Then
\begin{equation}\label{eq.tildeU}
    \maU(\maA\oplus \maL) \seq \mathbb{K} \oplus \tilde{\maU}(A\oplus L)\,.
\end{equation}

\begin{definition}
    \label{def.enveloping1}
    Let $\tilde{\maU}(\maA\oplus \maL)$ be as in Equation \eqref{eq.tildeU}.
    The \emph{enveloping algebra of the Lie-Rinehart algebra $(\maA, \maL)$} is,
    by definition, the quotient
    \begin{equation*}
        \maU(\maA, \maL) \ede \tilde{\maU}(\maA\oplus \maL) / I\,,
    \end{equation*}
    where $I$ is the ideal generated by elements of the form
    $i(b, 0)i(a, X) - i(ba, bX)$.
\end{definition}

The enveloping algebra $\maU(\maA, \maL)$ can be characterized by the following
universal property: if we are given an algebra $B$, an algebra homomorphism
$\phi_{\maA}: \maA \to B$, and a Lie algebra homomorphism
$\phi_{\maL}: \maL \to B$ (regarding $B$  as a Lie algebra
for the commutator) such that $\phi_{\maA}(a) \phi_{\maL}(X) = \phi_{\maL}(aX)$
and $[\phi_{\maL}(X),  \phi_{\maA}(a)]=  \phi_{\maA}(X(a))$,
then there exists a unique homomorphism $\Phi: \maU(\maA, \maL) \to B$ such
that $\Phi\circ i_{\maA}= \phi_{\maA}$ and $\Phi\circ i_{\maL}= \phi_{\maL}$.
This allows us to introduce the following important representation of
$\maU(\maA, \maL)$.

\begin{definition}[\cite{MoerdijkMrcunLieRine}]
    \label{def.morph.rho}
    Using the universal property of $\maU(\maA, \maL)$, we let
    \begin{equation*}
        \rho_{\maA, \maL}: \maU(\maA, \maL)\to \End(\maA)
    \end{equation*}
    be the unique morphism such that $\rho_{\maA, \maL}(i_\maA(a))$ acts by
    left multiplication by $a \in \maA$ and $\rho_{\maA, \maL}(i_\maL(\maG))$
    acts by derivations on $\maA$. We also let $\Diff(\maA, \maL) :=
    \rho_{\maA, \maL}(\maU(\maA, \maL))$ be the \emph{algebra of differential
    operators} associated to $(\maA, \maL)$.
\end{definition}

In case $(\maA, \maL) = (\maA(\maG), \maL(\maG))$ is associated
to the CF-groupoid $\maG$, then we shall write $\Diff(\maA, \maL) = \Diff(\maG)$.

In the examples considered in this article,
$\rho_{\maA, \maL}$ is injective and hence
\begin{equation*}
    \Diff(\maA, \maL) \ede
    \rho_{\maA, \maL} (\maU(\maA, \maL)) \simeq \maU(\maA, \maL)\,.
\end{equation*}
General injectivity results of this kind for $\rho_{\maA, \maL}$
will be included in \cite{BCNQ_LR}. If $\rho_{\maA, \maL}$ is injective,
we may regard  $ \Diff(\maA, \maL)$ as the subalgebra of $ \End(\maA)$
generated by $\maA$, acting as multiplication operators, and by $\maL$,
acting as derivations.
If, moreover, $(\maA, \maL)$ is associated to the CF-groupoid $\maG$, then from
Remark \ref{rem.inj.restr}, we have $\maA = \maA(\maG)_0$ acting injectively
on $\maA$ as multiplication operators, and $\maL =
\maL(\maG)_0\subset \maC(\units; \mathbf{A}(\maG))$ also acting injectively
as derivations, so $ \Diff(\maG)$ is isomorphic to the subalgebra of
$\End(\maA)$ generated
by functions in $ \maA(\maG)_0$ and vector fields in $ \maL(\maG)_0$, which
can then be regarded as differential operators of order $1$.

To obtain our final results on the Schr\"odinger operator, we will
need pseudodifferential operators. One important feature is that $\Diff(\maG)$
maps injectively to (it is even included in) the much larger class
$\Psi(\maG)$ of pseudodifferential operators
on $\maG$. Recall the following definition \cite{LMN2000, Paterson2000}.

\begin{definition}\label{def.pseudo}
    The set $\Psi^{m}(\maG)$ consists of
    right invariant, continuous families of pseudodifferential operators of order $m$ on
    $\maG_{x} := d^{-1}(x)$ with compactly supported convolution kernels
    $k_P$ (which are distributions on $\maG$).
\end{definition}

Similarly to the representation $\rho_{\maA, \maL}$, an operator
$P \in \Psi^{m}(\maG)$ acts on $\maA(\maG)$ via a map $\pi_{0}(P)$, and
we have $\Diff(\maG) := \Diff(\maA(\maG), \maL(\maG)) \subset
\pi_{0}(\Psi^{\infty}(\maG))$. The morphism $\pi_{0}$ is not injective in general,
but it will be in the examples considered in this paper.  
Note that we will not need
this class of operators up to section \ref{ssec:schroedinger}.

\subsection{CF-groupoids and associated Lie-Rinehart algebras over intervals}
\label{ssec.interval}
We now look at some of the simplest examples of Lie-Rinehart algebras.
The Lie-Rinehart algebras present some analytic subtleties, as it will be apparent
even from the one dimensional examples studied in this section. These examples will
also be an important building block in the constructions of the next section.
We agree that $\NN = \{1, 2, \ldots \}$ and $\ZZ_+ = \NN \cup \{0\}$. Let $\maC(K)$
denote the space of continuous functions on some compact space $K$, as before.

Let us introduce now some basic spaces. To this end, let
$N \in \ZZ_+ \cup \{\infty\}$ and
$\alpha, \beta \in [-\infty, \infty] =: \oR$,
$\alpha < \beta$. We consider
\begin{equation}\label{eq.def.CA}
    \begin{gathered}
        X(t) \ede \phi(t)\pa_t \,,\quad
        \phi : (\alpha, \beta) \to (0, \infty) \
        \mbox{ \emph{ continuous,}}\\
        \maC_{\phi}^{(N)} \ede  \{  u : [\alpha, \beta] \to \RR
        \mid  X^k u \in \maC([\alpha, \beta]) \; \text{ for all } \;
        0 \le k \le N \}\,,   \quad
        \mbox{and}\\
        \maV_{\phi}^{(N)} \ede \{ u X \mid u \in \maC_{\phi}^{(N)} \} \,.
    \end{gathered}
\end{equation}
We shall write $\maC_{\phi}^{(N)} = \maC_{\phi}^{(N)}([\alpha, \beta])$
and $\maV_{\phi}^{(N)} = \maV_{\phi}^{(N)}([\alpha, \beta])$
when we want to be more precise. The relation
``$X^k u \in \maC([\alpha, \beta])$'' is defined by induction to mean that
``$X^{k-1} u$ is differentiable on $(\alpha, \beta)$
and $X^ku := \phi (X^{k-1} u)'$ extends to a continuous function on $[\alpha, \beta]$''.
Note that $\maC_{\phi}^{(\infty)} = \cap_{n \in \ZZ_+} \maC_{\phi}^{(n)}$
and $\maC_{\phi}^{(0)} = \maC([\alpha, \beta])$.

We fix a function $\phi$ and a vector field $X$ as in Equation \eqref{eq.def.CA}
for the rest of this section. The main examples that will concern us will be
$\phi(t) = (t-\alpha)^a (\beta-t)^b$, with $a, b > 0$ (even $\ge 1$ for the most
part). Most of our example will be such that $\phi$ is differentiable on the open
interval $(\alpha, \beta)$, but this is by no means necessary for the development
of the theory, so it will not be assumed.
We have the following simple remark.

\begin{remark} \label{rem.easy}
    Let us assume that $\phi$ is smooth on $(\alpha, \beta)$.
    The assumption that $\phi(t) > 0$ gives then
    that every $u\in \maC_{\phi}^{(\infty)}$ is smooth on $(\alpha, \beta)$. In fact,
    we even have
    \begin{equation*}
        \CIc((\alpha, \beta)) \subset \maC_{\phi}^{(\infty)}
        \subset \maC([\alpha, \beta])
        \cap \CI((\alpha, \beta))\,.
    \end{equation*}
    If $\phi$ is smooth on a \emph{bounded} interval $[\alpha, \beta]$, then
    $\CI([\alpha, \beta]) \subset \maC_{\phi}^{(\infty)}$, but, in general,
    $\CI([\alpha, \beta]) \neq \maC_{\phi}^{(\infty)},$ and this is one of
    the point of considering Lie-Rinehart algebras: it allows for more
    general functions than the usual smooth functions.
\end{remark}

We have the following simple, technical result.

\begin{proposition}\label{prop.1D}
    We shall use the notation of Equation
    \eqref{eq.def.CA}. Let $n \in \ZZ_+ \cup \{\infty\}$. Then:
    \begin{enumerate}
        \item $\maC_{\phi}^{(n)}$ is an algebra.

        \item $\maV_{\phi}^{(n)}\maC_{\phi}^{(n+1)} \subset \maC_{\phi}^{(n)}$.

        \item $[\maV_{\phi}^{(n+1)}, \maV_{\phi}^{(n+1)}] \subset \maV_{\phi}^{(n)}$.
    \end{enumerate}

   In particular, $(\maC_{\phi}^{(\infty)},\maV_{\phi}^{(\infty)})$ is a
   Lie-Rinehart algebra.
\end{proposition}

\begin{proof} We proceed in order.  It suffices to assume $n\in \ZZ^+$, as the
    results for $n=\infty$ will then follow.
    (1) Let us show that $\maC_{\phi}^{(n)}$ is an algebra. This is clear if $n = 0$,
    so let us assume that $n > 0$. We have
    $X (fg) = X(f) g + f X(g) $, as $X := \phi \pa_t$ is a derivation. We obtain, for
    $k\in \ZZ_+$, $k \le n$,
    \begin{equation}\label{eq.Xk}
        X^k (uv) \seq \sum_{j=0}^n C_k^j
        \, X^j (u) X^{k-j}(v) \,,
    \end{equation}
    where $ C_k^j= \frac{k! }{j! (k-j)!} $, $0\leq j\leq n$ are the usual binomial
    coefficients. It follows that, if $X^j (u) $ and $X^{k-j}(v)$ are continuous on
    $[\alpha, \beta]$ for all $j = 0, \ldots, k$, then $X^k (uv)$ will also be
    continuous on $[\alpha, \beta]$. In particular, if
    $u, v\in \maC_{\phi}^{(n)}$ then $uv\in \maC_{\phi}^{(n)}$, so $\maC_{\phi}^{(n)}$
    is an algebra.

    (2) follows immediately from (1). Indeed, we have
    $X \maC_{\phi}^{(n+1)} \subset
    \maC_{\phi}^{(n)}$ by definition (Equation \ref{eq.def.CA}).
    Let $u \in \maC_{\phi}^{(n)}$
    and $v \in \maC_{\phi}^{(n+1)}$. Then
    \begin{equation*} (uX)(v) \ede u(X(v)) \in \maC_{\phi}^{(n)}\end{equation*}
    since $u, X(v) \in \maC_{\phi}^{(n)}$, which was proved to be
    an algebra.

    (3) also follows immediately from (1) and (2) in view of the
    Leibnitz rule, which yields the identity
    $[uX, vX] \seq [uX(v) - vX(u)] X\,.$

    In particular, $\maC_{\phi}^{(\infty)}$ is a commutative algebra,
    $\maV_{\phi}^{(\infty)}$ is a Lie algebra, by (3), which is also a left
    $\maA$-module, and acts by derivations on $\maC_{\phi}^{(\infty)}$, as
    $X$ does. Moreover, if $a\in \maC_{\phi}^{(\infty)}$ and
    $uX, vX \in \maV_{\phi}^{(\infty)}$, then
    \begin{equation*}
    [uX, avX]  \seq \big (uX(av) -a vX(u) \big ) X
    \seq (uX)(a)vX + a [uX, vX].
    \end{equation*}
    Hence, the Leibnitz rule is satisfied
    and $(\maC_{\phi}^{(\infty)},\maV_{\phi}^{(\infty)})$ is a
    Lie-Rinehart algebra.
\end{proof}

\subsection{Properties of the algebras $(\maC_{\phi}^{(\infty)},\maV_{\phi}^{(\infty)})$}
\label{ssec.more}
We continue to fix $\alpha, \beta \in \oR$ and a continuous
function $\phi : (\alpha, \beta) \to (0, \infty)$, as in the definition of
the Lie-Rinehart algebra $(\maC_\phi^{(\infty)}, \maV_\phi^{(\infty)})$ of
Equation \eqref{eq.def.CA}.

\begin{lemma}\label{lemma.same.LCIp}
    Let $\psi : (\alpha, \beta) \to (0, \infty)$ be a continuous function and 
    $\maC_{\psi}^{(\infty)}$ be defined as in \eqref{eq.def.CA}.
    \begin{enumerate}
        \item If $\phi/\psi \in \maC_{\psi}^{(\infty)}$, then
        $\maC_{\psi}^{(\infty)} \subset \maC_{\phi}^{(\infty)}.$
        \item Consequently, if $\phi/\psi \in \maC_{\psi}^{(\infty)}$ and
        $\psi/\phi \in \maC_{\phi}^{(\infty)}$, then
            $\maC_{\phi}^{(\infty)} = \maC_{\psi}^{(\infty)}.$
    \end{enumerate}
\end{lemma}

\begin{proof}
Let $Y := \psi \pa_t$. Then
$\maC_{\psi}^{(\infty)} := \{u : [\alpha, \beta] \to \CC \mid
Y^k u \in \maC_{\psi}^{(\infty)},\ k \in \ZZ_+\},$ by definition,
Equation \eqref{eq.def.CA}. Let
$X := \phi \pa_t = \frac{\phi}{\psi} Y.$ The assumption that $\frac{\phi}{\psi} \in
\maC_{\psi}^{(\infty)}$ gives that $X \in \maV_{\psi}^{(\infty)}$.
Since $\maV_{\psi}^{(\infty)}$ acts by derivations on
$\maV_{\psi}^{(\infty)}$, we have that $X^k u \in
\maV_{\psi}^{(\infty)} \subset \maC$, for
all $k \in \ZZ_+$, and hence $u \in \maV_{\phi}^{(\infty)}$, by
definition. The second statement follows from the first one, by symmetry.
\end{proof}

We shall need the following lemma.

\begin{lemma} \label{lemma.maybe.useful}
    Let $X := \phi \pa_t$, as before, and $n, j \in \ZZ_+$.
    \begin{enumerate}
        \item Let us assume $\phi' \in \maC_\phi^{(n+j-2)}$ if $n \ge 2$
        (we impose no condition if $n =0$ or $n = 1$). Then there
        exist $a_k \in \maC_{\phi}^{(j)}$, $k = 0, \ldots, n-1$, such that
        \begin{equation*}X^{n} u \seq \sum_{k=0}^{n-1} a_k \phi^{k} \pa_t^k u 
            + \phi^{n} \pa_t^n u\,.\end{equation*}

        \item Let us assume $\phi' \in \maC_\phi^{(n-2)}$ if $n \ge 2$.
        We have $u \in \maC_{\phi}^{(n)}$ if, and only if
        $\phi^{k}\pa_t^k u \in \maC([\alpha, \beta])$, for all $k = 0, 1, \ldots, n$.
    \end{enumerate}
\end{lemma}

\begin{proof}
    (1) Let $\maS_n := \{ a_j \neq 0 \}$ be the set of non-zero coefficients
    of $X^n = \sum_{k=0}^{n} a_k \phi^{k} \pa_t^k$.
    We shall prove the statement (i.e. $\maS_{n} \subset \maC_{\phi}^{(j)}$
    if $\phi' \in \maC_\phi^{(n+j-2)}$)
    by induction on $n$ (and for all $j$).
    We have $\maS_0 = \maS_1 = \{1\}$ by definition and $\maS_2 = \{1,\phi'\}$.
    Hence the statement is true for $n \le 2$ (for $n = 2$ this is because we have
    assumed that $\phi' \in \maC_{\phi}^{(j)}$).

    Let us prove the induction step. Let us assume $\phi' \in \maC_{\phi}^{(n+1+j-2)}$
    for some $n \ge 2$.
    Then $\maS_n \subset \maC_{\phi}^{(j+1)}$, by the induction hypothesis. Next,
    \begin{align*}
        X^{n+1} u & \seq X \Big( \sum_{k=0}^{n} a_k \phi^{k} \pa_t^k u \Big ) \\
        & \seq \sum_{k=0}^{n} \big [ X(a_k) + k a_k\phi' + a_{k-1} \big ] \phi^{k} \pa_t^k u
        + a_n \phi^{n+1} \pa_t^{n+1} u\,,
    \end{align*}
    (where we have set $a_{-1} := 0$).
    Hence $\maS_{n+1} \subset \maC_{\phi}^{(j)}$.

    (2) Follows by induction from definitions using (1) for $j = 0$.
\end{proof}

We obtain the following consequence.

\begin{corollary} \label{cor.description}
    Let $n \in \ZZ_+ \cup \{\infty\}$.
    We have $\phi' \in \maC_{\phi}^{(n)}$ if, and only if,
    $\phi^k \pa_t^{k+1} \phi \in \maC([\alpha, \beta])$ for all $0 \le k \le n$.
\end{corollary}

\begin{proof}
    The proof is again by induction on $n$. For $n = 0$ and $n = 1$ this is checked
    directly. Let us assume the statement to be true for some $n \ge 1$ and prove it
    for $n +1$. Assume $\phi' \in \maC_{\phi}^{(n+1)}$. Then Lemma
    \ref{lemma.maybe.useful} (2) for $u = \phi'$ gives
    $\phi^k \pa_t^k \phi' \in \maC([\alpha, \beta])$
    for all $0 \le k \le n+1$. Conversely, if $\phi^k \pa_t^k \phi' \in \maC([\alpha, \beta])$
    for all $0 \le k \le n+1$, then the induction hypothesis gives $\phi' \in \maC_{\phi}^{(n)}$
    Hence we can use Lemma \ref{lemma.maybe.useful} (2) again for $u = \phi'$ to conclude
    that, in fact, $\phi' \in \maC_{\phi}^{(n+1)}$.
\end{proof}

Similarly, we obtain the following consequence.

\begin{corollary} \label{cor.description2}
    Let us assume that $\phi^k \pa_t^{k+1} \phi \in \maC([\alpha, \beta])$
    for all $k \in \ZZ_+$
    and let $\eta \in \maC([\alpha, \beta])$ be such that
    $\phi^k \pa_t^k \eta \in \maC([\alpha, \beta])$
    for all $k \in \ZZ_+$. Then $\eta \in \maC_{\phi}^{(\infty)}$.
\end{corollary}

Let $\psi : [\alpha, \beta] \to [0, \infty)$, $\psi(t) > 0$
for $t \in (\alpha, \beta)$ be continuous.
Also, let
\begin{equation}\label{eq.defCphipsi}
C_{\psi, \phi}:= \frac{\phi \psi'}{\psi} = X(\psi)/\psi,
\end{equation}
which is defined
initially for $t \in (\alpha, \beta)$, but which we extend to $\alpha$ and $\beta$ by
continuity, when possible. We shall need yet the following corollary.

\begin{corollary} \label{cor.conseq}
    Let us assume that $\phi(t) = (t-\alpha)^{a}$ and $\psi(t) = (t-\alpha)^{b}$
    for $t$ close to $\alpha$ and that $\phi(t) = (\beta -t)^{a'}$ and
    $\psi(t) = (\beta - t)^{b'}$ for $t$ close to $\beta$.
    Assume $a, a' \ge 1$ and $b, b' \ge 0$. Then $  \psi \in\,
        \maC_{\phi}^{(\infty)}$ and also
            \begin{equation*}
        \phi'\, ,\   C_{\psi, \phi}  \ede \phi \psi'/\psi\, \
        %\psi
         \in\,
        \maC_{\phi}^{(\infty)}\,.
    \end{equation*}
\end{corollary}

\begin{proof}
The fact that $\psi \in \maC_{\phi}^{(\infty)}$
    follows directly from the definition since, for instance,
    near $\alpha$, $X^{k}\psi = b(b+a-1)\ldots (b+ (k-1)(a-1))(x-\alpha)^{
        b + k(a-1)},$ which is continuous.
    Recall that $C_{\psi, \phi} := \frac{\phi \psi'}{\psi} = X(\psi)/\psi$.
    Let us check the conditions of Corollaries
    \ref{cor.description} and \ref{cor.description2}
    near $\alpha$ if $b > 0$. Indeed, for $t \in [\alpha, \beta]$ close
    to $\alpha$, $C_{\psi, \phi}(t) = b(t-\alpha)^{a-1}$ and
    \begin{equation*}
        a^{-1}\phi^k \phi^{(k+1)}(t) \seq
        b^{-1}\phi^k \pa_t^{k} C_{\psi, \phi}(t) \seq
        (a-1)\ldots (a-k)(t-\alpha)^{(k+1)(a-1)}\,,
    \end{equation*}
    are both continuous near $\alpha$ (including $\alpha$).
    If $b = 0$, everything is the same, except the last
    formula, in which we have $C_{\psi, \phi} = 0$, and hence
    the same conclusion. The calculations near $\beta$ are similar.
\end{proof}

We remark that, if $\phi$ and $\psi$ are as in the last corollary, then
\begin{equation*} 
    \begin{cases}
        C_{\psi, \phi}(\alpha) \seq b & \mbox{ if } a \seq 1\\
        C_{\psi, \phi}(\alpha) \seq 0 & \mbox{ if } a > 1\,.
    \end{cases}
\end{equation*}
A similar relation holds true at $\beta$:
\begin{equation*} 
    \begin{cases}
        C_{\psi, \phi}(\beta) \seq b' & \mbox{ if } a' \seq 1\\
        C_{\psi, \phi}(\beta) \seq 0 & \mbox{ if } a' > 1\,.
    \end{cases}
\end{equation*}

As we have seen $(\maC_{\phi}^{(\infty)},\maV_{\phi}^{(\infty)})$ is a Lie-Rinehart algebra,
it is natural to ask when it comes from a CF-groupoid? We shall see that this is the case
if $\phi$ (or $X := \phi \pa_t$) is \emph{complete} in the sense that
\begin{equation}\label{eq.assumpt.phi}
    \int_{\alpha}^\gamma \frac1{\phi (t)}dt \seq \int_{\gamma}^{\beta}
    \frac1{\phi (t)}dt \seq \infty\
\end{equation}
for one (equivalently, for any) $\gamma \in (\alpha, \beta)$. We fix
from now on $\gamma \in (\alpha, \beta)$.

Let, for further reference,
\begin{equation}\label{eq.def.F}
    F(x) \ede \int_{\gamma}^{x} \frac{dt}{\phi(t)}\,.
\end{equation}
This notation will be fixed throughout the rest of the paper.

\begin{proposition} \label{prop.e.flow}
    Let us assume that our continuous
    function $\phi : (\alpha, \beta) \to (0, \infty)$
    satisfies the completeness conditions of Equation
    \eqref{eq.assumpt.phi}. Then $X := \phi \pa_t$
    generates a one parameter group of homeomorphisms $\sigma_s$,
    $s \in \RR$, of
    $[\alpha, \beta]$ that are $\maC^1$ on $(\alpha, \beta)$.
\end{proposition}

What this proposition says is that under the completeness condition
we have a well defined flow
\begin{equation}
    \sigma_s : [\alpha, \beta] \to [\alpha, \beta]\,,
    \quad \sigma_{s}\circ \sigma_{t} \seq \sigma_{s+t}\,, \ \ s, t \in \RR\,,
\end{equation}
which as usually means
\begin{equation}
    \frac{d}{ds} \sigma_s (x) \seq \phi(\sigma_s(x)) \,, \ \mbox{ for }
    x \in (\alpha, \beta)\,,
\end{equation}
and which is equivalent to the usual relation
\begin{equation}
    \frac{d}{ds} f(\sigma_s (x)) \seq [Xf](\sigma_s(x)) \,.
\end{equation}

\begin{proof}

Let $F$ be as in Equation \eqref{eq.def.F}.
Then, our assumptions on $\phi$ give that $F : (\alpha, \beta) \to (-\infty, \infty)$ is
$\maC^1$-homeomorphism. Then the desired flow $\sigma_s$ is given by
\begin{equation}
    \sigma_s(x) \seq F^{-1}(F(x) + s)\,.
\end{equation}
Note that even though $\phi$ is not Lipschitz, the flow $\sigma_s$ is
well defined. Indeed, Any two solutions $t(s)$ and $\tilde t(s)$ of the
differential equation
\begin{equation*}
    y'(t) \seq \phi(y(t)),
\end{equation*}
with $y(0)=\tilde y (0) = \gamma$ must satisfy $F(y(s))=F(\tilde y(s))$,
for all $s\in \mathbb{R}$. By assumption we have $\phi(x)>0$, for
$x\in(\alpha,\beta)$. For this reason, $F$ is a strictly monotone
increasing function and therefore one-to-one. Hence the two solutions
$y$ and $\tilde y$ must agree and we get uniqueness of the Cauchy problem
and a well-defined flow on the open interval $(\alpha,\beta)$.

Since $\lim_{x \to \alpha} F(x) = -\infty$ and
$\lim_{y \to -\infty } F^{-1}(y) = \alpha$,
we obtain that $\lim_{x \to \alpha} \sigma_s(x) = \alpha$. Similarly,
$\lim_{x \to \beta} \sigma_s(x) = \beta$. We shall then set
\begin{equation}\label{eq.extension}
    \sigma_s(\alpha) \ede \alpha \ \mbox{ and }\
    \sigma_s(\beta) \ede \beta\,.
\end{equation}
In this way, the resulting maps $\sigma_s$ define a one-parameter group of
\emph{homeomorphisms} of $[\alpha, \beta]$ that are $\maC^1$-homeomorphisms
in the interior $(\alpha, \beta)$.
\end{proof}

Let us give some concrete examples, i.e., start with the simplest example,
namely, the one that underpins the so called ``$b$-calculus'' of Melrose
\cite{MelroseActa} and Schulze \cite{SchulzeBook91}. See also
\cite{GrieserBCalc, LauterSeiler, LeschBCalc, SchroheFrechet}.

\begin{example} \label{example.b.calc}
    Let $(\alpha, \beta) = (0, \infty)$ and $\phi(t) = t$. Then the
    function $F$ of Equation \eqref{eq.def.F} is given by $F(x) = \ln x$,
    which is a diffeomorphism $F : (0, \infty) \to \RR$ with
    $F^{-1}(y) =  \exp(y) = e^y $. Hence the flow $\sigma_s$, $s \in \RR$,
    of Proposition \ref{prop.e.flow} will be given by
    \begin{equation*}
        \sigma_s(x) \seq F^{-1}(F(x) + s)
        \seq \exp(\ln x + s) \seq e^s x\,.
    \end{equation*}
    We note that, in this example, $\sigma_s$ \emph{is smooth also at} $0$.
\end{example}

This example can be generalized as follows.

\begin{example} \label{example.ca.calc}
    Let $\alpha = 0$, $\beta > \gamma = 1$, and $a > 1$ (the case $a = 1$
    was treated in the previous example and here $a$ is \emph{not}
    assumed to be an integer!). Let $\phi(t) = t^a$ for $t \le 1$ and arbitrary
    for $t > 1$, as long as it satisfies the needed conditions,
    namely $\phi$ is continuous and $> 0$ on $(\alpha, \beta)$
    (Equation \ref{eq.def.CA}) and Equation \eqref{eq.assumpt.phi}. Then
    \begin{equation*}
        \begin{cases}
            F(x) \seq \frac{ 1-x^{{1-a}}}{a-1} & \mbox{ for } x \le 1 \mbox{ and}\\
            F^{-1}(y) \seq [1-(a-1) y]^{\frac1{1-a}} & \mbox{ for } y \le 0\,.
        \end{cases}
    \end{equation*}
    Therefore, if $x \le 1$ and $F(x) + s \le 0$, we obtain
    \begin{align*}
        \sigma_s(x) & \seq F^{-1}(F(x) + s)\\
        & \seq [1 - (a-1)s x^{a-1}]^{\frac1{1-a}} x\,.
    \end{align*}
    In particular, $\sigma_s : [0, \beta] \to [0, \beta]$
    is $\maC^\infty$ at $0$ if, and only if, $a$ is an integer.
\end{example}

Yet another generalization of Example \ref{example.b.calc} is to a \emph{bounded}
interval.

\begin{example} \label{example.b.calc2}
    Let $[\alpha, \beta] = [-1, 1]$ and $\phi(t) = (1+t)(1-t)$.
    Then $F(x) = \frac12 \ln \left| \frac{1+x}{1-x} \right|$
    is a diffeomorphism $F : (-1, 1) \to \RR$ with $F^{-1}(y) = \tanh(y)$.
    Therefore
    $\sigma_s(x) = \tanh( F(x) + s)$ and, in this case, $\sigma_s$ extends to a
    $\maC^\infty$ diffeomorphism $[-1, 1]\to [-1, 1]$, for all $s$.
\end{example}

The point is that if $\phi$ satisfies the completeness condition of Equation
\eqref{eq.assumpt.phi}, we can use the flow $\sigma_s$ to construct a CF-groupoid
whose Lie-Rinehart algebra is exactly
$(\maC_{\phi}^{(\infty)},\maV_{\phi}^{(\infty)})$. More precisely, consider
the group-action groupoid:
\begin{equation}\label{eq.def.Gphi}
    \maG_\phi \ede [\alpha, \beta] \rtimes_{\sigma} \RR\,,
\end{equation}
with the usual structural maps (see \cite{ConnesFol, ConnesNCG, RenaultB,
WilliamsBook})
\begin{equation}\label{eq.def.str.Gphi}
    \begin{cases}
        \ d(x, t) \ede x & \\
        \ r(x, t) \ede \sigma_t(x) &\\
        \ (y, s)(x, t) = (x, s + t) & \mbox{ if } \sigma_t(x) = y\,,
    \end{cases}
\end{equation}
where $\sigma_s$ is the flow given by Proposition \ref{prop.e.flow}.
The particular case $\phi \equiv 1$ and $[\alpha, \beta] = [-\infty, + \infty] =
\oR$ will play a special role; the associated groupoid (namely $\maG_\phi$) will
also be denoted $\maG_\phi = \maG_1$, for simplicity. So $\maG_1 := \oR \rtimes \RR$,
with the action of $\RR$ on $\oR$ by translations.

\begin{theorem} \label{theorem.descr.A}
    Let us assume that $\phi : (\alpha, \beta) \to (0, \infty)$ is
    continuous and that it satisfies the completeness condition of Equation
    \eqref{eq.assumpt.phi}, as before. Then the group-action groupoid $\maG_\phi$
    of Equations \eqref{eq.def.Gphi} and \eqref{eq.def.str.Gphi} is a
    CF-groupoid such that $F : [\alpha, \beta] \to \oR$ induces an
    isomorphism \begin{equation*}\maG_\phi \simeq \maG_1\,.\end{equation*} We hence have
    $\maA(\maG_\phi)_{0} = \maC_{\phi}^{(\infty)}$ and
    $\maL(\maG_\phi)_{0}= \maV_\phi^{(\infty)}$ (see Equation \eqref{eq.def.CA}
    and Remark \eqref{rem.new.r} for notation). Let us assume that the interval
    $[\alpha, \beta]$ is bounded and it is endowed with the usual smooth
    structure, then $\maG_\phi$ is a Lie groupoid if, and only if, the function
    $\phi$ is smooth on $[\alpha, \beta]$.
\end{theorem}

\begin{proof}
    That $\maG_\phi$ is a CF-groupoid follows directly 
    from~\cite[Proposition 5]{Paterson2000} (see the Example \ref{examples}\eqref{item.cor} 
    for the needed statement). 

    The relations $\maA(\maG_\phi)_{0} = \maC_{\phi}^{(\infty)}$ and
    $\maL(\maG_\phi)_{0}= \maV_\phi^{(\infty)}$ can be checked directly
    if $[\alpha, \beta]= \oR$ and $\phi \equiv 1$. They are hence true in general
    since the diffeomorphism $F : (\alpha, \beta) \to \RR$ of Equation \eqref{eq.def.F}
    maps $X := \phi \pa_t$ to $\pa_t$, by construction,  hence
    $\maC_{\phi}^{(\infty)}= \maC_1^{(\infty)}$,
    $\maV_\phi^{(\infty)}= \maV_1^{(\infty)}$. Moreover,
    $\maG_\phi \simeq \maG_1\,$, so the associated Lie-Rinehart algebras
    are also isomorphic.

    If the groupoid $\maG_\phi$ is a Lie groupoid, then $r$ must be smooth, and hence
    the action $\sigma_t$ (with $\sigma_t$ the flow of Proposition \ref{prop.e.flow})
    will be smooth on $[\alpha, \beta]$. Consequently, $\phi$ will also be smooth on
    $[\alpha, \beta]$. Conversely, if $\phi$ is smooth on $[\alpha, \beta]$, then the
    standard results on the differentiability of solutions of ODEs gives that $\sigma_t$
    is smooth on $[\alpha, \beta]$.
\end{proof}

\begin{corollary}\label{cor.if.ab}
    For example, assume in addition that $a, b \ge 1$ are
    such that $\phi(t) = (t-\alpha)^a$ for $t$ close to $\alpha$ and $\phi(t) = (\beta-t)^b$
    for $t$ close to  $\beta$. Then $\maG_\phi$ is a Lie groupoid if, and only if,
    $a, b \in \NN$.
\end{corollary}

Our results raise the following integrability question:
``Which Lie-Rinehart algebras are associated to some CF-groupoids?''
The gluing construction of \cite{NistorJapan}
(integrating on each orbit and then gluing) seems to be still relevant
here. In the following, we shall provide some constructions of CF-groupoids
associated to some given Lie-Rinehart algebras particularly useful for
studying Schr\"odinger operators with singular potentials.

\section{The basic fibered boundary case}

\label{sec:general_integration+Schroedinger}
We now construct the CF-groupoids used to study Schr\"odinger operators.
Our construction relies on the one dimensional examples
of the previous section. It is an instance of ``integrating''
a suitable Lie-Rinehart algebra.

\subsection{Formulation of the problem and main integration results}
\label{ssec.bs.fp}

We now use the results proved so far to investigate more complicated Lie-Rinehart
algebras and their associated CF-groupoids. Let $M_1$ be a smooth, compact manifold of
dimension $d$ (i.e. a closed manifold) and, for the rest of this section, we let
\begin{equation*}
    M \ede [\alpha, \beta] \times M_1 \,.
\end{equation*}
We continue to consider a continuous function $\phi : (\alpha, \beta) \to
(0, \infty)$ that is subject to the same assumptions as before
(namely, Equations \ref{eq.def.CA} and \ref{eq.assumpt.phi}). To fix  ideas,
we shall also assume that $\phi$ is \emph{smooth on $(\alpha, \beta)$
and extends by continuity to $[\alpha, \beta]$ such that $\phi(\alpha) = \phi(\beta)
= 0$.} We then  lift $\phi$ to a function $[\alpha, \beta] \times M_1 \to [0, \infty)$,
still denoted $\phi$, by abuse of notation. In addition to the function $\phi$, we consider
a second function $\psi$ with properties similar to those of $\phi$ (and those of the
function $\psi$ in the previous subsection), namely,
\begin{equation}\label{eq.def.psi}
    \begin{cases}
        \ \psi \in \maC_{\phi}^{(\infty)}\\
        \ \psi(t) > 0 \mbox{ on } (\alpha, \beta)\\
        \ C_{\phi, \psi}
        \ede \phi(t)\psi'(t)/\psi(t) \
        \mbox{ extends to a function in } \maC_{\phi}^{(\infty)}\,.
    \end{cases}
\end{equation}
(Recall that $\maC_{\phi}^{(\infty)}$ was introduced in Equation \eqref{eq.def.CA}
and that $\maC_{\phi}^{(\infty)} \subset \maC([\alpha, \beta])$.)
For simplicity, we shall also assume that $\alpha, \beta \in \RR$
(that is, that $[\alpha, \beta]$ is bounded).

Let us begin with a heuristical formulation of our problem.
Let us place ourselves in a local coordinate chart, so
$M_1$ is an open subset of $\RR^d$. We shall use the functions
$\phi$ and $\psi$ described in the previous paragraph.
We consider vector fields of the form
\begin{equation*}
    X \ede Y_0 \ede \phi \pa_t \,, \
    \mbox{ and }\ Y_k \ede \psi \pa_{x_k}\,, \quad 1 \le k \le d\,,
\end{equation*}
on $[\alpha , \beta] \times M_{1} \ni (t, x)$
and \emph{we will construct a CF-groupoid whose Lie-Rinehart algebra
is generated by these vector fields.}

To formulate more precisely our problem, let us introduce the following
notation. Let $\pi_0 : [\alpha, \beta] \times M_1  \to [\alpha, \beta]$ and
$\pi_1 : [\alpha, \beta] \times M_1  \to M_1$ be the
two projections. Then we have a canonical isomorphism of vector bundles
\begin{equation}\label{eq.direct.sum.tang}
    TM \simeq \pi_0^*(T[\alpha, \beta]) \oplus \pi_1^*(TM_1) \simeq
    \RR \times \pi_1^*(TM_1)\,,
\end{equation}
where the direct sum $\oplus$ is the direct sum of vector bundles on $M$ and $\times$
is the usual Cartesian product of sets. A section of $\pi_0^*(T[\alpha, \beta])$ will be
called a \emph{horizontal vector field} (on $M$) and a section of $\pi_1^*(TM_1)$ will be
called a \emph{vertical vector field} (on $M$). Accordingly, a vector field $Z$ on $M$
will decompose as
\begin{equation} \label{eq.def.u01}
    Z \ede (Z_0, Z_1) \in \pi_0^*(T[\alpha, \beta]) \oplus \pi_1^*(TM_1)\,,
\end{equation}
where $Z_0$ is a vector field on $M := [\alpha, \beta] \times M_1$ in the direction
of $[\alpha, \beta]$ and $Z_1$ is a vector field on $M$ in the direction of $M_1$.
Thus $Z_0$ is a horizontal vector field and $Z_1$ is a vertical vector field.
The copy of $\RR$ on the right hand side is due to the fact that every
horizontal vector field  $Z_0$ is of the form $v_0 \pa_t$, for some suitable function
$v_0$ on $M$, so, at each point of $M$, it is determined by a single scalar (in $\RR$).
We shall often use in what follows the isomorphisms
$\maA(\maG) \simeq \maA(\maG)_0 \subset \maC(\units)$ and $\maL(\maG) \simeq
\maL(\maG)_0 \subset \maC(\units, \mathbf{A}(\maG))$ defined
by restriction to units (as in Remark \eqref{rem.new.r}).

\begin{problem}
    \label{problem.integration}
    To find a CF-groupoid $\maG_{\psi, \phi}$ with units
    $\maG_{\psi, \phi}^{(0)} = M := [\alpha, \beta] \times M_1$
    whose associated Lie-Rinehart algebra of vector fields $\maL(\maG_{\psi, \phi})_0$
    is given by restriction to units (see Remark \eqref{rem.new.r}) by vector fields of
    the form
    \begin{equation*}
        Z \seq (\phi Z_{0}, \psi Z_{1})\,,
    \end{equation*}
    with $Z_{0}$ horizontal and $Z_{1}$ vertical vector fields on $M$
    (see Equation \eqref{eq.def.u01} for notation; additional smoothness properties may be
    required of $Z_{0}$ and $Z_{1}$, in particular, $\phi$ will be as in
    Equation \eqref{eq.def.CA} and $\psi$ will be as in Equation \eqref{eq.def.psi}).
\end{problem}

Of course, our Problem \ref{problem.integration} is an ``integration problem''
for Lie-Rinehart algebras. More precisely, let $(\maA(\maS), \maL(\maS))$
be the Lie-Rinehart algebra generated by vector fields of the form
$Z \seq (\phi Z_{0}, \psi Z_{1})$, as above. Then we will construct a CF-groupoid
$\maG_{\psi, \phi}$ integrating $(\maA(\maS), \maL(\maS))$, that is, such
that $(\maA(\maS), \maL(\maS)) = (\maA(\maG), \maL(\maG))$.

Our construction of the groupoid $\maG_{\psi, \phi}$
depends on whether or not $\psi$ vanishes at the end points ($\alpha$ and $\beta$).
The result, however, is independent of this condition. More precisely,
\emph{we will construct two CF-groupoids $\maH_{\psi}$ and $\maG_{\psi, \phi}$
with units $M$ satisfying, in turn, he following two results}
(recall the notation and isomorphisms of Remark \eqref{rem.new.r}).
The main properties of the groupoid $\maH_{\psi}$ are described by
the following proposition.

\begin{proposition} \label{prop.descr.A}
    Let $\psi$ be as in Equation \eqref{eq.def.psi}. There exists
    an explicit CF-groupoid $\maH_{\psi}$ (given by Equation \eqref{eq.def.Hpsi1},
    if $\psi(\alpha) > 0$, and by Equation \eqref{eq.def.Hpsi2}, if $\psi(\alpha)=0$),
    which is minimal (and hence unique up to isomorphism)
    such that the Lie-Rinehart algebra of $\maH_{\psi}$ is given by
    \begin{equation*} %\label{eq:descr.Afirst}
        \begin{gathered}
            \maA(\maH_{\psi})_0 \seq \{u \mid \psi^{j} Y_1 \ldots Y_j u \in \maC(M)\,,
            \ j \in \ZZ_+ \, \mbox{ and }\, Y_i \in \CI(M; \pi_1^*(TM_1)) \}\,\\
            %\label{eq:descr.Lfirst}
            \maL(\maH_{\psi})_0 \seq \psi \maA(\maH_{\psi}) \CI(M; \pi_1^*(TM_1))
            \subset \maC(M; TM) \,.
        \end{gathered}
    \end{equation*}

\end{proposition}

In turn, the main properties of the groupoid $\maG_{\psi, \phi}$ are described by
the following proposition. Before stating it, let us notice that the vector fields
$Y_i$ above are \emph{vertical} vector fields on $M$.
In particular, the CF-groupoid $\maH_{\psi}$ integrates the vertical vector fields
on $M$ that are multiples of $\psi$. In order to simplify the statement
of the result satisfied by $\maG_{\psi, \phi}$, it is useful
to notice that, in the last proposition, we can replace $\CI(M; \pi_1^*(TM_1))$
(the set of smooth, vertical vector fields on $M$) with $\CI(M_1; TM_1) \subset
\CI(M; \pi_1^*(TM_1))$ (the smooth, vertical vector fields on $M =
[\alpha, \beta] \times M_{1}$ that are independent of $t \in [\alpha, \beta]$).

\begin{theorem} \label{thm.descr.A}
    Let $\phi$ and $\psi$ be as in Equations \eqref{eq.def.CA}
    and \eqref{eq.def.psi} and $\maH_{\phi}$ be the CF-groupoid
    of Proposition \ref{prop.descr.A}. Then there exists an action
    of $\RR$ on $\maH_{\psi}$ whose infinitesimal generator is
    $X = \phi \pa_{t}$ and the action (or semi-direct product) groupoid
    $\maG_{\psi, \phi} :=
    \maH_{\psi} \rtimes \RR$ is a CF-groupoid with Lie-Rinehart algebra given by
    \begin{equation*}%\label{eq:descr.Asecond}
        \begin{gathered}
            \maA(\maG_{\psi, \phi})_0 \seq \{u
            \mid \psi^{j} X^k Y_1 \ldots Y_j u \in \maC(M)\,,
            \ j, k \in \ZZ_+ \, \mbox{ and }\,
            Y_i \in \CI(M_1; TM_1) \}\\
            %\label{eq:descr.Lsecond}
            \maL(\maG_{\psi, \phi})_0 \seq \phi \maA(\maG_{\psi, \phi})\pa_t
            + \psi \maA(\maG_{\psi, \phi}) \CI(M_1; TM_1) \subset \maC(M; TM) \,.
        \end{gathered}
    \end{equation*}
    If $\frac{\psi}{\psi_{1}}, \frac{\psi_{1}}{\psi} \in \maC_{\phi}^{(\infty)}$,
    then we have canonical isomorphisms $\maH_{\psi} \simeq \maH_{\psi_{1}}$
    and $\maG_{\psi, \phi} \simeq \maG_{\psi_{1}, \phi}$.
\end{theorem}

Let $\psi, \psi_{1} : [\alpha, \beta] \to [0, \infty)$.
In view of the above theorem, we shall write
\begin{equation}\label{eq.def.sim.phi}
    \psi \sim_{\phi} \psi_{1}\ \Leftrightarrow\ \frac{\psi}{\psi_{1}},
    \frac{\psi_{1}}{\psi} \in \maC_{\phi}^{(\infty)}\,.
\end{equation}
(More precisely, the quotients $\frac{\psi}{\psi_{1}}$ and
$\frac{\psi_{1}}{\psi}$ extend by continuity to functions defined on
$[\alpha, \beta]$ and that these extensions are in $\maC_{\phi}^{(\infty)}$.)

\begin{remark}\label{rem.generation}
    It follows that the Lie-Rinehart algebra $(\maA(\maG_{\psi, \phi}),
    \maL(\maG_{\psi, \phi}))$ associated to the CF-groupoid $\maG_{\psi, \phi}$
    is ``generated'' in an obvious sense by $X = \phi \pa_{t}$
    and $\psi \CI(M_1; TM_1)$.
\end{remark}

The following two subsections are devoted to proving the above two
results assuming $\psi(\beta) > 0$. The proof will be split in two
according to the values of $\psi$ at $\alpha$. (The case $\psi(\beta) = 0$
is completely similar and can also be obtained from the case
$\psi(\beta) > 0$ using a gluing procedure.)

\subsection{The case $\psi(\alpha) > 0$}
We thus consider first the case when our given continuous function
$\psi: [\alpha, \beta] \to [0, \infty)$ is such that
$\psi(\alpha) \psi(\beta) \neq 0$. Then $\psi > 0$ on $[\alpha, \beta]$.
In this case, our construction will independent of $\psi$ with these
properties. Our first step is to construct a groupoid $\maH_{\psi}$
of Proposition \ref{prop.descr.A}. Explicitly,
\begin{equation}\label{eq.def.Hpsi1}
    \maH_{\psi} \ede [\alpha, \beta] \times (M_1 \times M_1)\,,
\end{equation}
with structural maps:
\begin{equation}\label{eq.def.dr0}
    \begin{gathered}
        d(s, x, y) \ede (s, y) \quad \mbox{and}
        \quad r(s, x, y) \ede (x, s) \,,
        \quad \mbox{and}\\
        (s, x, y) (s, y, z) \ede (s, x, z)
        \qquad s \in [\alpha, \beta]\, \mbox{ and }\, x, y \in M_{1}\,.
    \end{gathered}
\end{equation}
In the standard terminology (see, for instance, \cite{CNQ}), $\maH_{\psi}$
is a Lie groupoid, product of the space $[\alpha, \beta]$ (a groupoid that
has only units) with the pair groupoid $M_1 \times M_1$.
(The pair groupoid $M_{1} \times M_{1}$ is obtained by removing the factor
$s$ in the above formulas, see Example \ref{examples}\eqref{item.pair}).
The minimality of $\maH_{\psi}$ follows from the fact that the pair groupoid
$M_{1} \times M_{1}$ is the minimal groupoid with algebroid $TM_{1}$.
Recall that $M := [\alpha, \beta] \times M_{1}$
and a vector field on $M$ is \emph{vertical} if it is tangent to all the
copies $\{c\} \times M_{1}$ of $M_{1}$, $c \in [\alpha, \beta]$.
We see right away that $\maH_{\phi}$ is independent on $\psi$.

We want to prove Proposition \ref{prop.descr.A} in this case. Since $\psi$ is
non-vanishing and continuous the required identifications of Proposition
\ref{prop.descr.A} become
\begin{equation}
    \label{eq:descr.A1first}
        \maA(\maH_{\psi})_0 \seq \{u \mid \psi^{j} Y_1 \ldots Y_j u \in \maC(M)\}
        \seq \Cio(M)
\end{equation}
(for all $j \in \ZZ_+$ and  $Y_i \in \CI(M; \pi_1^*(TM_1))$) and
\begin{equation}
    \label{eq:descr.L1first}
        \maL(\maH_{\psi})_0 \seq \psi \maA(\maH_{\psi}) \CI(M; \pi_1^*(TM_1)) \seq
        \psi \Cio(M; \pi_1^*(TM_1)) \,.
\end{equation}
This amounts to the fact that $\maA(\maH_{\psi})_0$ consists of continuous
functions on $[\alpha, \beta ] \times M_1$ that are smooth in the direction of $M_{1}$
(the vertical direction) and $\maL(\maH_{\psi})_0$ consists of continuous
vertical vector fields on $M$ that are smooth in the direction
of $M_1$ (again the vertical direction).

\begin{proof}[Proof of Proposition \ref{prop.descr.A} for $\psi(\alpha)\psi(\beta) \neq 0$]
    We notice that $\maH_{\psi}$, $\maA(\maH_{\psi})_0$,
    $\maL(\maH_{\psi})_0$, and the other spaces of Equations
    \eqref{eq:descr.A1first} and \eqref{eq:descr.L1first}
    ($\Cio(M)$ and $\psi \Cio(M; \pi_1^*(TM_1))$) are independent
    on $\psi > 0$ because $1/\psi \in \maC_{\phi}^{(\infty)}$ (which is due, in turn,
    to the fact that $\psi \in \maC_{\phi}^{(\infty)}$ and $\psi > 0$ on $[\alpha, \beta]$).
    It is enough then to assume that $\psi = 1$.
    By definition, $\maA(\maH_{\psi})_0$ consists of longitudinally smooth
    right-invariant functions. For this Lie groupoid we have $d(s, x, y) = (s, y)$,
    where $s \in [\alpha, \beta]$ and $x, y \in M_{1}$. Hence, given a function
    $f(s, x, y)$, it must be smooth in the $x$-variable and continuous in the other
    ones. Now if $f$ is right-invariant, then
    \begin{equation*}
        f(s, x, y) \seq f((s, x, y)(s, y,z)) \seq f(s, x,z)\,.
    \end{equation*}
    From here, we see that right-invariant functions are actually independent of
    the second variable, and thus we can identify them with functions $f(s, x)$ on
    $M$, which are smooth along the variable $x \in M_1$. This is exactly the first
    set given in Proposition~\ref{prop.descr.A}.

    Let us now identify the second set in Proposition~\ref{prop.descr.A}. First, we
    have that a vector field $X \in \maL(\maH_{\psi})_0$ must be a derivation along
    the $d$-fibers. In local coordinates, this means that it must be of the form:
    \begin{equation*}
        X \seq \sum_{i=1}^{\dim M_1} X_i(s, x, y)\pa_{x^i} \,.
    \end{equation*}
    where $X$ a priori is only smooth in the $x$ variables. Now we need to check
    what conditions on $X_i(s, x, y)$ are imposed by the right invariance. For this,
    we need to compute the differential $(R_{(s, y, z)})_*$ of the right multiplication.
    Since, $(s, x, y)(s, y, z) = (s, x, z),$
    we have
    \begin{equation*}
        (R_{(s, y, z)})_* =
        \begin{pmatrix}
            id & 0 & 0 \\
            0 & 1 & 0 \\
            0 & 0 &  0
        \end{pmatrix}.
    \end{equation*}
    Thus
    \begin{equation*}
        (R_{(s, y, z)})_* X(s, x, y) = \sum_{i=1}^{\dim M_1} X_i(s, x, z)\pa_{x^i}.
    \end{equation*}
    If $(R_{(s, y,z)})_* X(s, x, y)  = X(s, x, y)$, then we must have
    $X_i(s, x, y) = X_i(s, x,z)$,
    which again allows us to conclude that $X_{i}$ is independent of
    the last variable, and hence that $X_i \in \maA(\maH_{\psi})_0$.
    This completes the proof.
\end{proof}

Let $\sigma$ be the action of $\RR$
on $[\alpha, \beta]$ associated to $\phi$ as in the previous section
(See Proposition \ref{prop.e.flow} and the discussion following it). We
extend this action to an action of $\RR$ on
$\maH_{\psi} := [\alpha, \beta] \times M_1  \times M_1$
that is trivial on $M_1 \times M_1$. We then let
\begin{equation}\label{eq.def.maG.psi0}
    \maG_{\psi, \phi} \ede \maH_{\psi} \rtimes_{\tau} \RR
    \simeq M_1 \times M_1 \times \maG_{\phi} \,,
\end{equation}
(see Theorem \ref{thm.descr.A}).
We see from the definitions that $\maG_{\psi, \phi} = \maG_{1, \phi}$ is
independent of $\psi$, as long as $\psi$ is smooth and it does not vanish
at the end points.

\begin{proof}[Proof of Theorem \ref{thm.descr.A} for $\psi(\alpha)\psi(\beta) \neq 0$]
The groupoid $\maG_{\psi, \phi} := \maG_{\psi} \rtimes \RR$ of Theorem \ref{thm.descr.A}
(see also Equation \eqref{eq.def.maG.psi0}) is the Cartesian product
of the action groupoid $\maG_{\phi}$ and the pair groupoid $M_1 \times M_1$.
We shall identify it, as a set, with $M_1 \times M_1 \times [\alpha, \beta] \times \RR$ via
\begin{equation}\label{eq.def.maG.psi1}
    \maG_{\psi, \phi} \ede \maH_{\psi} \rtimes_{\tau} \RR
    \simeq M_1 \times M_1 \times \maG_{\phi} \simeq
    M_1 \times M_1 \times [\alpha, \beta] \times \RR\,.
\end{equation}
that is, the product of the pair groupoid $M_1 \times M_1$ with the
action groupoid $\maG_\phi$ introduced in the previous section.
See Equations \eqref{eq.def.str.Gphi}, \eqref{eq.def.str.Gphi}, and \eqref{eq.def.dr0}
for the products on these groupoids. Hence, we have
\begin{equation*}
    (x, y,  s,  t)(y, z, \sigma_{-\tau}(s), \tau) \seq
    (x, z, \sigma_{-\tau}(s), t + \tau) \,.
\end{equation*}
In particular, if a function $f$ is right-invariant, then
\begin{equation*}
    f(x, y, s, t) = f(x, z, \sigma_{-\tau}(s), t + \tau).
\end{equation*}
As a consequence if $f(x, y, s, t)$ is right invariant, then,
\begin{itemize}
    \item it is independent of $y$ and,
    \item it is entirely determined by the restriction to $t + \tau = 0$.
\end{itemize}
Thus $f(x, y, s, t) = f_{0}(x, \sigma_{t}(s))$, for some function
$f_{0}$ on the set of units $M_1 \times [\alpha, \beta] \simeq M$
(the diffeomorphism is given by switching $[\alpha, \beta]$ and $M_{1}$).
Since $f\in \maA(\maG_{\psi, \phi})_0$ must be smooth along the $d$-fibres,
it should be differentiable in the $x$ and $t$ variables. To establish the
isomorphism from the statement, we notice that the reduced function
$f_{0}(x,\sigma_{t}(s))$ depends only on $r(x, y, s, t)$. Hence it can be
seen as the pull-back $r^*u$ of a function $u\in \CI (M)$. Fix $t$ and
consider a new variable $t'(s) = \sigma_{t}(s)$ and extend it by continuity
to $t'(-\infty) = \alpha$, $t'(-\infty) = \beta$. Then
$f(x,\sigma_{t}(s)) = f(x,t')$ and the differentiability
in $t$ transfers to differentiability via the vector field
$\phi(t') \partial_{t'}$.
The differentiability properties now follow.

Next we prove the second isomorphism. The vector fields in
$\mathcal{L}(\mathcal G_{\psi,\phi})_0$ must be tangent to the
$d$-fibers and right invariant. Since $d(x,y,\tau,t) = (y,\tau)$,
those vector fields must be derivations in $t$ and $x$. Vector
fields $\partial_t$ and $V\in \CI(M_1,TM_1)\subset \CI(M,TM)$
can be seen tautologically as vector fields on the $d$-fibers. Note that
all of them are already right-invariant: $\partial_t$ because the right
multiplication is just a translation in the $t$-variable, and $V$ by an
argument similar to the proof of Proposition~\ref{prop.descr.A}
for our case.
Since those vector fields generate any other $d$-tangent vector field,
all the invariant vector fields must be linear combinations of those
ones with invariant functions as coefficients. So it only remains to take
$t'=\sigma_{t}(s)$ as a new variable and we get exactly
$\maA(\maG_{\psi, \phi})_0$. In the case $\psi > 0$, all the groupoids
and associated spaces are independent on $\psi$. This gives
\begin{equation}
    \maG_{\psi_{1}, \phi} \seq \maG_{1, \phi} \seq \maG_{\psi, \phi}\,,
    \qquad \psi, \psi_{1} > 0\,,
\end{equation}
that is, the isomorphism of the last statement in Theorem \ref{thm.descr.A}.
\end{proof}

\subsection{The case $\psi(\alpha) = 0$} Recall that $\psi(\beta) > 0$.
\label{ssec.Hpsi}
As before, we begin by defining the
CF-groupoid $\maH_{\psi}$ integrating $\psi \CI(M; \pi_1^*(TM_1))$. Let
\begin{equation} \label{eq.def.Hpsi2}
    \maH_{\psi} \ede TM_1 \times \{\alpha\} \cup M_1 \times M_1
    \times (\alpha, \beta]
\end{equation}
as a set. The groupoid structure is the following. The domain, range, and multiplication maps
on  $M_1 \times M_1 \times (\alpha, \beta]$ are induced from the inclusion
\begin{equation*}
    M_1 \times M_1 \times (\alpha, \beta] \times \RR \subset
    \maG_{1, \phi} \ede M_1 \times M_1 \times [\alpha, \beta] \times \RR
\end{equation*}
in the groupoid $\maG_{1, \phi}$ of the previous subsection. Let $\pi_{TM_1} : TM_1 \to M_1$
be the canonical projection. It will be convenient to write
$M = [\alpha , \beta] \times M_{1}$ (we have switched the two factors).
On the other hand, if $s =\alpha$, these maps are given by
\begin{equation}\label{eq.def.dr}
    \begin{cases}
        d(v, \alpha) \seq r(v, \alpha) \ede (\pi_{TM_1}(v), \alpha)
        \in M_1 \times \{\alpha\} \subset M & \mbox{ and}\\
        (v, \alpha) (w, \alpha) \ede (v + w, \alpha) \,,&
    \end{cases}
\end{equation}
where $v, w \in TM_1$ are vectors in the same fiber over $M_1$, that is,
$\pi_{TM_1}(v) = \pi_{TM_1}(w)$. Thus, $\maH_{\psi}$ is a ``family'' of groupoids indexed
by $s \in [\alpha, \beta]$, such that, to $s \in (\alpha, \beta]$ we associate the pair
groupoid $M_1 \times M_1$ and to $s = \alpha$ we associate the tangent bundle
$TM_1$ with the associated groupoid structure of a family of Lie groups.

To define the topology and CF-groupoid structure on $\maH_{\psi}$, we consider
a variant of the ``deformation to the normal cone'' considered by Connes, Debord, Skandalis,
and others \cite{ConnesNCG, DebordSkandalis}. Thus, we first consider on
$TM_1 \times \{\alpha\}$ and on $M_1 \times M_1 \times (\alpha, \beta]$ the natural smooth
(manifold) structures and the Lie groupoid structures. To glue these smooth
structures, we use the exponential coordinates of \cite{NistorJapan} to define
the smooth structure in the neighborhood of a point $(v, \alpha) \in TM_1 \times \{\alpha\}$.
This is done as follows. Then, let us fix some arbitrary connection on $M_1$ (that is,
on its tangent space $TM_1$) and let $\exp : U \subset TM_1 \to M_1$ be its associated
exponential map \cite{Petersen} (here $U$ is its maximal domain
of definition, a subset of $TM_1$). Let $\pi_{TM_1} : TM_1 \to M_1$ be
the projection map. It is known that
\begin{equation}\label{eq.def.chart.Psi}
    \Psi \ede (\pi_{TM_1}, \exp) : U \to M_1 \times M_1
\end{equation}
is a diffeomorphism from an open neighborhood of the zero section of $TM_1$
onto an open neighborhood of the diagonal $\delta_{M_1} \subset M_1
\times M_1$. Let $v \in TM_1$. Let us define now a neighborhood of $(v, \alpha)$
in $\maH_{\psi}$. Let first $1 \ge s_{0} > 0$ be such that $s v \in U$ for all
$|s| \le s_{0}$. Then, let us choose a relatively compact neighborhood $W$ of $sv$
in $TM_1$ such that $[-1, 1]\overline{W} \subset U$. (So the exponential $\exp(sw)$
is defined for all $w \in W$ and $|s| \le 1$.)
Let $\epsilon > 0$ be such that $\psi(s) < s_{0}$ if $\alpha \le s < \alpha + \epsilon < \beta$.
We then let $\Phi : W \times [0, \epsilon) \to \maH_{\psi, \phi}$,
\begin{equation}\label{eq.coord.ch.Phi}
    \Phi(w, s) \ede
    \begin{cases}
        \ (w, \alpha)  & \ \mbox{ if }\ s = 0 \\
        \ (\Psi ( \psi(s) w), \alpha + s) \in M_1 \times M_1 \times (\alpha, \beta)
        & \ \mbox{ if }\ s > 0\,.
    \end{cases}
\end{equation}
(We recall from Equation \eqref{eq.def.chart.Psi}
that $(\Psi ( \psi(s) w), s') = (\pi_{TM_1}(w), \exp( \psi(s)w), s')$.)
We stress the appearance of the factor $\psi(s)$ in these formulas. Except this
``rescaling term'', our formulas are as in the aforementioned papers.
Finally, the desired neighborhood of $v$ consists of
\begin{equation}\label{eq.def.neigh.V}
    V \ede \Phi(W \times [0, \epsilon))
\end{equation}
with the induced topological and continuous family of manifolds structure
from $\Phi$. (The continuous family of manifolds structure on $V$ comes from
the map $d$.)
We are ready now to prove Proposition \ref{prop.descr.A} for
$\psi(\alpha) = 0$ and $\psi(\beta) > 0$.

\begin{proof}[Proof of Proposition \ref{prop.descr.A} for $\psi(\alpha)\psi(\beta) = 0$]
Similar to the proof of this Proposition in the case $\psi(\alpha) > 0$,
the invariant functions for $s \in (\alpha, \beta)$ satisfy
\begin{equation*}
    u(x, y, s) = u_0(x, s)
\end{equation*}
and must be smooth in $x$ for all $s \in (\alpha, \beta)$.
The difference now is in the behavior as $s$ approaches $\alpha$ or $\beta$.
We shall concentrate on the case $s$ close to $\alpha$, the case
$s$ close to $\beta$ being completely similar. So let us change
the coordinates using the map $\Phi$. In this coordinates the domain and
the range maps are given by
\begin{align*}
    d(w,t) &\seq \pi_{TM_1}(w)\\
    r(w,t) &\seq \exp(\psi(t)w)\,,
\end{align*}
where $(w, t) \in W \times [0, \epsilon) \subset TM_1 \times [0, \epsilon)$
(see the definition of the map $\Phi$, Equation \eqref{eq.coord.ch.Phi})
So differentiating along $d$-fibres in these coordinates is differentiating
along vertical fibres of $w\in TM_1$. Since $x = \exp(\psi(t)w)$ we have that
any derivation in the $w$ variables will be a derivation in $x$ pushforwarded
by the differential of $r$, that is, by
\begin{equation}\label{eq:differential}
    d(\exp (\psi(t)\cdot)) \seq \psi(t) d \exp\, .
\end{equation}
Note that $d$ here indicates the differential. Since the exponential map
is a diffeomorphism for sufficiently small $t$,
we have that the derivative along the $d$-fiber exists if, and only if, $\psi Y u$
is continuous for any derivation $Y\in \CI(M_1,TM_1)$. The claimed
description of $\maA(\maH_{\psi})_0$ is then obtained by noticing
that $\psi$ and $Y$ commute (i.e. $Y(\psi) = 0$).

Now we characterize $\maL(\maH_\psi)_0$. Again as in the proof of this
Proposition in the case $\psi(\alpha) > 0$,
the right invariant vector fields for $t>\alpha$ agree with continuous families of
derivations in the $x$ variables parametrized by the $t$ variable. The only thing to
understand is what happens as we approach $t=\alpha$. For this we change to the $w$-variables
using the diffeomorphism $\Phi$ of Equation \eqref{eq.coord.ch.Phi}.
If $X_t \in \CI(M_1,TM_1)$ is a family of such derivations, then in the $w$ variables we
get via~\eqref{eq:differential}:
\begin{equation*}
(\exp(\psi(t))\cdot)^{-1}_* X_t = \frac{1}{\psi(t)}(\exp )^{-1}_*X_t.
\end{equation*}
Hence, we must have that $X_t/\psi(t)$ is continuous in $t$. From here the
characterization follows.
\end{proof}

We want to define an action of $\RR$ on $\maH_{\psi}$ integrating $X = \phi \pa_t$
to define $\maG_{\psi, \phi} := \maH_{\psi} \rtimes \RR$ Theorem \ref{thm.descr.A}.
First, the automorphisms $\sigma_s$, $s \in \RR$, extend (as for the case $\psi > 0$
in the previous subsection)
to $M_1 \times M_1 \times (\alpha, \beta]$. Recall the spaces
$\maC_{\phi}^{(\infty)} \subset \maC([\alpha, \beta])$
introduced in Equation \eqref{eq.def.CA}.

\begin{proposition} \label{prop.extension}
    We continue to assume that $\psi(\alpha) = 0$ (as throughout this
    subsection). Then the assumption $C_{\psi, \phi} \ede  \frac{X(\psi)}{\psi}
    \ede \frac{\phi \psi'}{\psi} \in \maC_{\phi}^{(\infty)}\,,$
    of Equation \eqref{eq.def.psi} (see also Equation \eqref{eq.def.CA})
    allows us to conclude that $\sigma_s$ extend to CF-groupoid automorphism
    of $\maH_{\psi}$, for all $s \in \RR$. This extension to the end point
    $\alpha$ is such that
    \begin{equation*}
        \sigma_s(v, \alpha) \seq (e^{-C_{\psi, \phi}(\alpha) s} v, t)\,.
    \end{equation*}
\end{proposition}

\begin{proof}
    Let $s \in \RR$. Let us check the continuity and longitudinal
    smoothness of $\sigma_s$ near $\alpha$ (near the other points
    $(\alpha, \beta]$, it is obvious; it follows, for instance from
    case already proved of Theorem \ref{thm.descr.A}). Let
    $\lambda := C_{\psi, \phi}(\alpha)$. Let $\Phi$
    and $\Psi$ be as in Equations \eqref{eq.def.chart.Psi} and
    \eqref{eq.coord.ch.Phi}. We have
    \begin{equation}
        \sigma_s(\Phi(w,t)) \seq
        \begin{cases}
            \ (e^{-\lambda s} w, \alpha)  & \ \mbox{ if }\ t = \alpha \\
            \ (\Psi ( \psi(t) w), \sigma_s(t)) \in M_1 \times M_1 \times (\alpha, \beta)
            & \ \mbox{ if }\ t \in (\alpha , \alpha + \epsilon)\,.
        \end{cases}
    \end{equation}
    Hence, the map $\sigma_s$ maps the fibers of $d : \maH_{\psi} \to M$
    to themselves smoothly. In other words, $\sigma_s$ is longitudinally
    smooth. To complete the proof, we need to show that it is continuous.
    We have
    \begin{equation*}
        \Phi(\psi(t)w,\sigma_s(t))  \seq
        \left(\pi_{TM_1}(w),\exp \left(\frac{\psi(t)}{\psi(\sigma_{s}(t))}
        \psi(\sigma_{s}(t)) w\right),\sigma_s(t)\right)
    \end{equation*}
    So it only remains to show that $\lim_{t \to \alpha}
    \frac{\psi(t)}{\psi(\sigma_s(t))} = e^{-\lambda s}$. Indeed, we can compute:
    \begin{align*}
        \ln \frac{\psi(\sigma_s(t)) }{\psi(t)} & \seq \int_{0}^s \pa_\tau
        \ln \psi(\sigma_\tau(t))
        d\tau
        \seq \int_{0}^s \frac{\pa_\tau \psi(\sigma_\tau(t))}{\psi(\sigma_\tau(t))}
        d\tau \\
        & \seq \int_{0}^s \frac{(X\psi)(\sigma_\tau(t))}{\psi(\sigma_\tau(t))} d\tau
        \seq \int_{0}^s C_{\psi, \phi} (\sigma_\tau(t)) d\tau \to
        C_{\psi, \phi}(\alpha)s
        \seq \lambda s
    \end{align*}
    because $\sigma_\tau(t) \to \alpha$ as $t \to \alpha$ uniformly in
    $\tau \in [0, s]$.
\end{proof}

\begin{proof}{Proof of Theorem \ref{thm.descr.A} for $\psi(\alpha) = 0$
    and $\psi(\beta) > 0$.} The proof is very similar to that for the
    case $\psi(\alpha) > 0$. (All the difficulties were transfered to
    the definition of $\maH_{\psi}$.)
    The action of $\RR$ on $\maH_\psi$, for $t \in (\alpha, \beta]$, agrees
    exactly with the action of $\RR$ on the groupoid~\eqref{eq.def.maG.psi1}.
    From the proof of the case $\psi(\alpha)$ positive
    of Proposition~\ref{prop.descr.A} it follows, that the
    invariant functions for $t>\alpha$ can be identified with functions of the form
    $f(x,y,t,s)=f(x,\sigma_s(t))$
    (notice that $t$ and $s$ have been switched compared to the proof of
    Proposition~\ref{prop.descr.A}).
    As in the proof of Theorem~\ref{thm.descr.A} for $\psi(\alpha) = 0$ by compactifying the
    $\RR-$action and introducing a new variable $t'=\sigma_s(t)$, the
    differentiability in $s$, implies that $\phi_{t'}\partial_{t'}f$ must be continuous.
    On the other hand, the differentiability in $x$ is exactly as in
    Proposition~\ref{prop.descr.A}. The same holds for higher derivatives and
    hence the characterization follows.

    The characterization of $\maV_{\psi,\phi}(\maG_{\psi,\phi})_0$ is again very similar
    to the proofs of the above propositions. We have that elements of this space must
    be differentiations in $(x,s)$. For derivations in $x$ we repeat word-by-word the
    argument in the proof of Proposition~\ref{prop.descr.A} and find that
    $\psi \maC^{\infty}(M_1,TM_1)$ belong to $\maV_{\psi,\phi}(\maG_{\psi,\phi})_0$. Next we
    consider a derivation in $s$. It is a differentiation of a smooth family of
    automorphisms of $\maH_\psi$ by Proposition.~\ref{prop.extension}. Hence it is enough
    to describe this differentiation for $t \in (\alpha, \beta]$ and
    then extend it to whole of $\maH_\psi$. The argument in the proof of
    Theorem~\ref{thm.descr.A} when $\psi(\alpha)\neq 0$ says that for $t\in (\alpha, \beta]$ 
    the derivative
    $\partial_s$ is pushforwarded to $\phi(t')\partial_{t'}$ after a change of variables.
    Any other invariant vector field must be a linear combination of mentioned
    ones with invariant functions as coefficients.

    Finally, again, as a set $\maG_{\psi, \phi}$ is independent of $\psi$
    and the identity map $\maG_{\psi, \phi} \to \maG_{\psi_{1}, \phi}$ is an
    isomorphism of (plain) groupoids. If $\psi \sim_{\phi} \psi_{1}$
    (see Equation \eqref{eq.def.sim.phi} for the notation), then the identity
    map $\maG_{\psi, \phi} \to \maG_{\psi_{1}, \phi}$ is also a $\Cio$-map. This
    completes the proof.
\end{proof}

\begin{remark}
   Theorem \ref{theorem.descr.A} can be cast in a more general framework as
    follows. Let $\maH$ be a CF-groupoid and let $G$ be a Lie group acting by
    $\Cio$-homeomorphisms on $\maH$ (thus they are continuous and smooth on each
    of the fibers of the maps $d$ or $r$). Then $\maG \ede \maH \rtimes G$
    is a CF groupoid (see also Example \ref{examples} \eqref{item.cor} and \eqref{item.thm}).
    We have that $\maA(\maG)_{0}$ identifies with the smooth
    vectors in $\maA(\maH)$ and $\maL(\maG)_{0}$ identifies with the
    direct sum of the set of smooth vectors in $\maL(\maH)_{0}$ with
    $\maA(\maG)_{0} Lie(G)$.
\end{remark}

The case of
a general interval $[\alpha, \beta] \subset \overline{\RR}$ is dealt
with by using a homeomorphism of $[\alpha, \beta]$ to a bounded
interval, provided that this homeomorphism is a diffeomorphism
in the interior.

As already mentioned, the remaining case $\psi(\beta) = 0$
can be treated similarly. Better, yet, it can be treated by localization
and by gluing groupoids, as in \cite{CarvalhoComeQiao, Come}.
We recall now the gluing construction, for the
benefit of the reader.

\subsection{Localization and gluing of CF-groupoids}
\label{ssec.gluing}
In this section, the symbol $\maG$ (possibly decorated with various indices)
will denote a CF-groupoid. The set of units of $\maG_i$ will be denoted
$\maG_i^{(0)}$, unless otherwise stated. The domain and range maps will be
denoted by $d$ and, respectively, $r$, as usual.
%
%\begin{proposition} \label{prop.localization}
    Let $\maG$ be a CF-groupoid with units $\units$ and
    $U \subset \units$ be an open subset. Let
    \begin{equation*}
        \maG_U^U \ede d^{-1}(U) \cap r^{-1}(U)\,,
    \end{equation*}
    Then $\maG_U^U$ is a CF-groupoid with units $U$.
%\end{proposition}
We can glue CF-groupoids in a similar way as it was done for Lie groupoids in
\cite{CarvalhoComeQiao, Come}.

\begin{theorem} \label{thm.gluing}
    Let $M$ be a locally compact space and $\maG_i$, $i = 1, \ldots, N$,
    be CF-groupoids with units open subsets
    $\maG^{(0)}_i \subset M$ such that $\cup_{i=1}^N \maG^{(0)}_i =M$.
    We assume that these groupoids satisfy the gluing conditions:
    \begin{enumerate}
        \item (Isomorphism on the common part)
        There are CF-groupoid isomorphisms
        \begin{equation*}
            f_{ji}:(\maG_j)_{\maG^{(0)}_i \cap \maG^{(0)}_j}^{\maG^{(0)}_i
            \cap \maG^{(0)}_j}
             \rightarrow
            (\maG_j)_{\maG^{(0)}_i \cap \maG^{(0)}_j}^{\maG^{(0)}_i
            \cap \maG^{(0)}_j},
        \end{equation*}
        such that $f_{ij}=f^{-1}_{ji}$ and $f_{ij}f_{jk}=f_{ik}$ on the common domains,
        for all $i,j,k\in \{1, 2, \cdots, N\}$,.
        \item (The orbits are included completely in at least one groupoid) For each
        $x\in M$, there exist $i_x \in \{1, 2, \cdots, N\}$, such that
        \begin{equation*}
            \maG_j\cdot x:=\{r(g) \,| \, d(g)=x \} \subset \maG^{(0)}_{i_x}
        \end{equation*}
        for all $j \in \{1, 2, \cdots, N\}$.
    \end{enumerate}
    Then their union $\maG := \bigsqcup_{i=1}^N \maG_i/\sim $, where $\sim$ is the
    equivalent relation generated by $g \sim f_{ji}(g)$, for all
    $i, j\in \{1, 2, \cdots, N\}$ and $g\in \maG_i$, is a CF-groupoid.
    The Lie-Rinehart algebra of $\maG$ can be identified as $\maA(\maG)_{0}
    \seq \{u \in \maC(M) \mid u\vert_{\maG^{(0)}_i} \in \maA(\maG_{i})_{0}\}$,
    and  $\maL(\maG)_{0} \seq \{u \in \maC(M, \mathbf{A}(\maG)) \mid
    u\vert_{\maG^{(0)}_i} \in \maL(\maG_{i})_{0}\}$.
\end{theorem}

The proof goes along the same lines as the one in \cite{CarvalhoComeQiao}.
Details will be included in \cite{BCNQ_LR}.
To use the gluing theorem, we observe that the reduction of the
groupoid $\maG_{\psi, \phi}$ to any relatively compact subset of $(\alpha, \beta)$
is independent of $\psi$ and $\phi$. More precisely, if $I := [\alpha', \beta']
\subset (\alpha, \beta)$, then
\begin{equation*}
    (\maG_{\psi, \phi})_{I}^{I} \simeq (I \times M_{1})^{2}
\end{equation*}
(the pair groupoid, see \ref{examples}\eqref{item.pair}). Note that this relies on
the minimality of $\maH_{\psi}$.

\subsection{The groupoid associated to Schr\"odinger operators}

Let us consider now $M = S^{n-1} \times [0, \infty]$ (thus $[\alpha, \beta] = [0, \infty]$)
and four parameters $a, b, a', b'$, with $a, a' \ge 1$ and $b, b' \ge 0$. We then
choose our functions $\phi$ and $\psi$ to satisfy
\begin{equation}\label{eq.cond.phi.psi}
    \begin{cases}
        \phi(t) \seq t^a & \mbox{ if } t \in  [0, \epsilon/2]\\
        \phi(t) \seq t^{-a'} & \mbox{ if }
        t \in  [2\epsilon^{-1}, \infty]\\
        \psi(t) \seq t^b & \mbox{ if } t \in  [0, \epsilon/2]\\
        \psi(t) \seq t^{-b'} & \mbox{ if }
        t \in  [2\epsilon^{-1}, \infty]\,,
    \end{cases}
\end{equation}
for some $0 < \epsilon < 1$. (This is, of course, in addition to the usual assumptions:
$\phi$ is strictly positive and continuous on $(\alpha, \beta)$ and $\psi$ satisfies
the conditions \eqref{eq.def.psi} on $[\alpha, \beta]$.) Moreover, for simplicity, we
shall assume that $\phi$ is smooth on $(\alpha, \beta)$. Up to an isomorphism, the
resulting groupoid $\maG_{\psi, \phi}$ is independent of $\epsilon$ and of the choices of
$\phi$ and $\psi$ satisfying Equation \eqref{eq.cond.phi.psi}, by the last part
of Theorem \ref{thm.descr.A}. Therefore, it makes sense to set
\begin{equation}\label{eq.def.maSab}
    \maS \seq \maS_{a,a',b,b'} \ede \maG_{\psi, \phi}\,,\qquad
    \phi, \psi \mbox{ satisfying Equation \eqref{eq.cond.phi.psi}\,.}
\end{equation}

We now look at some specific properties of the groupoid $\maS = \maS_{a,a',b,b'}$.
Recall that, for a CF-groupoid $\maG$, the algebra of differential
operators $\Diff(\maG)$ was defined as the image of the enveloping algebra of
$\Diff( \maA(\maG)_{0}, \maL(\maG)_{0})$ in $\End(\maA(\maG)_{0})$, see
Subsection \ref{ssec:universal_envelop}.
Let
    \begin{equation}\label{eq.def.2rho}
    \begin{cases}
        \rho_{0}, \rho_{\infty} : S^{n-1} \times [0, \infty] \to [0, 1] &\\
        \rho_{0}(x, t) = t \ \mbox{ and } \ \rho_{\infty}(x, t) = 1\  & \mbox{ for }
        t < \epsilon\\
        \rho_{\infty}(x, t) = 1/t \ \mbox{ and } \ \rho_{0}(x, t) = 1
        & \mbox{ for }
        t > 1/\epsilon\,,
    \end{cases}
\end{equation}
be smooth, depend only on the second variable, and be such that
$\rho_{0}$ vanishes only when $t = 0$
and $\rho_{\infty}$ vanishes only when $t = \infty$.
Thus these two functions are the defining functions of the two
components of the boundary of $S^{n-1} \times [0, \infty]$. In
view of the notation \eqref{eq.def.sim.phi}, we can thus take
\begin{equation}\label{eq.sim.with.rho}
    \phi \seq \rho_{0}^{a}\rho_{\infty}^{a'} \quad \mbox{ and }
    \quad \psi \seq \rho_{0}^{b}\rho_{\infty}^{b'} \,.
\end{equation}
(The function $\phi$ used in the introduction is the same as
the function $\phi$ used here.)

\begin{proposition}\label{proposition.diff.ops}
    We shall denote generically by $D_j \seq Y_{1} Y_{2} \ldots Y_{j}$ a
    differential monomial of order $j$ in smooth, vertical vector fields on 
    $M$ (so $Y_{i}$ are tangent to all
    $S^{n-1} \times \{ t \}$) and $\maS = \maS_{a, a', b, b'}$, 
    with  $\phi \seq \rho_{0}^{a}\rho_{\infty}^{a'}$ and 
    $\psi \seq \rho_{0}^{b}\rho_{\infty}^{b'}$. Then
    \begin{equation}\label{cond}
        u \in \maA(\maS)_{0} \ \ \Leftrightarrow \ \
        i, j \in \ZZ_{+}\,,\ \ \phi^{i} \psi^{j}
        \,
        \pa_{t}^{i} D_j u \in \maC(M)\,,
    \end{equation}
    Similarly, a differential
    operator on the interior of $M$ is in $\Diff(\maS)$ if, and only if,
    it is a linear combinations of terms of the form
    $c_{i, j} \phi^{i} \psi^{j}\, \pa_{t}^{i} D_j$,
    where $c_{i, j} \in \maA(\maS)_{0}$, $i, j \in \ZZ_{+}$.
\end{proposition}

\begin{proof}
    This follows from Theorem \ref{thm.descr.A}, more precisely, from the characterization
    of the Lie-Rinehart algebras in that theorem, and from Corollary \ref{cor.description2}.
    (See also Remark \ref{rem.generation}.)
\end{proof}

\begin{definition} \label{def.cab}
    Let $\maS = \maS_{a, a', b, b'}$ and let $D$ be a differential operator on
    $M := S^{n-1} \times [0, \infty]$, as above. We shall say that $D$ is in the
    $c_{a,b}$-calculus (respectively, in the $c_{a',b'}$-calculus) if it coincides
    with a differential operator in $\Diff(\maS)$ in a neighborhood of $M_{1} \times \{0\}$
    (respectively, in a neighborhood of $M_{1} \times \{\infty\}$).
\end{definition}

Related calculi to the $c_{a,b}$ calculi were used in see, for instance in
\cite{Kottke-Lin, 
VasyReg, VasySurv}, also in
connection to Schr\"{o}dinger operators. See also \cite{DLR, 
kottkeInvent,
Mazzeo91, MelroseScattering}.

Let $d$ and $r$ be the structural maps of the groupoid
$\maS = \maS_{a, b, a', b'}$. We shall need the following property of
this groupoid.

\begin{proposition}\label{prop.quotient}
    Let $\rho_{0}$ and $\rho_{\infty}$ be the functions
    introduced in Equation \eqref{eq.def.2rho}, then
      \begin{equation}\label{eq.def.zeta}
        \zeta_{0} \ede \frac{\rho_{0}\circ d}{\rho_{0} \circ r}
        \quad \text{and} \quad \zeta_{\infty} \ede
        \frac{\rho_{\infty}\circ d}{\rho_{\infty} \circ r}
    \end{equation}
    extend to nowhere vanishing $\Cio$ functions on $\maS_{a, b, a', b'}$.
\end{proposition}

\begin{proof}
    The defining formula $\maG_{\psi, \phi} := \maH_{\psi} \rtimes_{\sigma} \RR$
    (Theorem \ref{thm.descr.A}), we see that the quotients $\zeta_{0}$
    and $\zeta_{\infty}$ depend only on $[0, \infty] \rtimes \RR$,
    via the projection $\maH_{\psi} \to [0, \infty]$.
    The required smoothness then follows from the explicit formulas
    for $\sigma_{s}(x)/x$ of Examples \ref{example.b.calc} (for $a = 1$
    or $a' = 1$) or \ref{example.ca.calc} (for $a > 1$
    or $a' > 1$).
\end{proof}

Consider the class $\Psi(\maG)$ of pseudodifferential operators on $\maG$,
Definition \ref{def.pseudo}. The previous lemma then yields the following result.

\begin{proposition}\label{prop.conj.rho}
    Let $\rho_{0}$ and $\rho_{\infty}$ be the functions
    introduced Equation \eqref{eq.def.2rho}, $\maS = \maS_{a, b, a', b'}$,
    $P \in \Psi^{m}(\maS)$. Let also $\zeta_{0}$ and $\zeta_{\infty}$ be as in
    Proposition \ref{prop.quotient} and $t, t' \in \RR$. Then the operator
    \begin{equation*}
        \tilde P := \rho_{\infty}^{-t'}\rho_{0}^{-t}
        P \rho_{0}^{t} \rho_{\infty}^{t'}
    \end{equation*}
    has distribution kernel $k_{\tilde P} = \zeta_{0}^{t} \zeta_{\infty}^{t'}k_{P}$,
    and hence it belongs to $\Psi^{m}(\maS)$ as well.
\end{proposition}

\begin{proof}
    The proof is a direct consequence of definitions.
    Indeed, the distribution kernel of $fPg$ is
    $k_{fPg}(\gamma) = f(r(\gamma))k_{P}(\gamma) g(d(\gamma))$.
    This gives the claimed form of the convolution kernel of
    $\tilde P$. Since $\zeta_{0}$ and $\zeta_{\infty}$ are
    $\Cio$ and $> 0$, by Proposition \ref{prop.quotient},
    $\zeta_{0}^{t} \zeta_{\infty}^{t'}$
    is $\Cio$ on $\maS$. This gives the result by the identification
    of the convolution kernels $k_{P}$ in \cite{LMN2000}.
\end{proof}

\section{An application to Schr\"odinger operators}
\label{ssec:schroedinger}
We now apply the results of the previous sections to study
Schr\"{o}dinger operators with power-law singularities at 0
and at $\infty$ introduced in Equation \eqref{eq.schrV}.
We continue to consider the functions $\rho_{0}$ and
$\rho_{\infty}$ of Equation \eqref{eq.def.2rho} (introduced right
before Proposition \ref{prop.quotient}). Recall that
the operator
\begin{equation*}
    H \ede -\Delta + V
\end{equation*}
of Equation \eqref{eq.schrV} is such that the potential
$V$ is of the form $V = \rho_{0}^{-2\gamma} \rho_{\infty}^{2\gamma'}V_0$, with $V_0$
smooth on $S^{n-1} \times [0, \infty]$ in spherical coordinates
$(\rho, x')$, $\rho \in [0, \infty]$, $x' \in S^{n-1}$. (We can relax
the smoothness assumptions on $V_{0}$, but the most general assumptions on
$V_{0}$ are much more difficult to state.)

We begin by studying the operator $-\Delta + V$ in a neighbourhood of
the origin, where $\rho_{0} = \rho$. We may also assume $\rho_{\infty} = 1$
there. Let us assume first that $\gamma \le 1$.
Proposition \ref{proposition.diff.ops} and the last part of
Corollary \ref{cor.conseq} show that the operator
\begin{equation}\label{eq.schr3}
    \rho^{2}(-\Delta + V) \seq (\rho \pa_\rho)^2
    + (n-2) \rho \pa_\rho + \Delta_{S^{n-1}} +  \rho^{2-2\gamma}V_0
\end{equation}
is in the $c_{1,0}$-calculus of the previous section (see Definition \ref{def.cab});
this is the case even if $2\gamma$ is not an integer due
to Corollary \ref{cor.conseq}. (That is, our coefficients \emph{do not}
have to be smooth up to the boundary. By contrast, $\rho^{2}(-\Delta + V)$
is in the $b$-calculus precisely if $2\gamma \in \ZZ$ and $2 \gamma \le 2$.)
On the other hand, if $\gamma \ge 1$ (but
keeping the assumption that $V_0$ is smooth in polar coordinates), we
have that
\begin{equation}\label{eq.schr4}
    \rho^{2\gamma} (- \Delta + V) \seq (\rho^{\gamma} \pa_\rho)^2
    + (n-1-\gamma) \rho^{\gamma-1} \rho^\gamma \pa_\rho +
    \rho^{2\gamma-2}\Delta_{S^{n-1}} + V_0
\end{equation}
is in the $c_{\gamma, \gamma-1}$ calculus, since its coefficients, $\rho^{\gamma-1}$
and $\rho^{2\gamma-2}$, are in the algebra $\maC_{\rho^\gamma}^{(\infty)} =
\maC_{\rho^\gamma}^{(\infty)} ([0, \infty))$
of Equation \eqref{eq.def.CA} by Corollary \ref{cor.conseq}.

We can carry a similar analysis at $+\infty$. Our assumption is that our
potential $V$ is of the form
$V(\rho x') \seq \rho^{2\gamma'}V_{0}'(\rho^{-1}, x')$,
with $V_{0}'$ smooth on $[0, \infty) \times S^{n-1}$.
Let us set $r := \rho^{-1}$ and let us write the Schr\"odinger operator
in the coordinates $(r, x') \in [0, \infty) \times S^{n-1}$
(that is, ``near $\infty$''). We have
$\rho \partial_\rho = - r\partial_r$. This gives us
\begin{equation}\label{eq.schr5}
    - \Delta + V \seq (r^2\pa_r)^2 -(n-1) r (r^2\pa_r)
    + r^2\Delta_{S^{n-1}}+r^{-2\gamma'}V_0'\,.
\end{equation}

Similar to the analysis at $\rho = 0$ we see that if $\gamma' \leq 0$,
then at the infinity this operator is in the $c_{2,1}$-calculus of Definition
\ref{def.cab}. If $\gamma' \ge 0$, then we factor out $r^{-2\gamma'}$. It is
a straightforward computation to check that:
\begin{equation*}
    (r^{2+\gamma'}\pa_r)^2 \seq r^{2\gamma'}(r^2\pa_r)^2 +
    \gamma' r^{2\gamma + 3}\pa_r\,.
\end{equation*}
Then
\begin{equation}\label{eq.schr6}
    r^{2\gamma'}(- \Delta + V) \seq (r^{2+\gamma'} \pa_r)^2
    - (n-1+\gamma') r^{1+\gamma'} (r^{2+\gamma'} \pa_r) +
    r^{2\gamma' + 2}\Delta_{S^{n-1}} + V_0'
\end{equation}
is in the $c_{2+\gamma',\gamma' + 1}$-calculus. We have thus
proved the following proposition.

\begin{proposition}\label{prop.inPsiS}
    Let $\rho_{0}$ and $\rho_{\infty}$ be the functions
    introduced right before Proposition \ref{prop.quotient},
    in particular,  $\rho_{0}(t)=t$ for $t$ small and
    $\rho_{\infty}(t)=t^{-1}$ for $t$ large, see
    Equation \eqref{eq.def.2rho}. Let $\gamma, \gamma' \in \RR$ and $V$
    be such that $\rho_{0}^{-2\gamma}\rho_{\infty}^{-2\gamma'}V$
    be smooth on $[0, \infty] \times S^{n-1}$. Let $\tilde \gamma :=
    \max\{\gamma, 1\}$ and, similarly, $\tilde \gamma' :=
    \max\{\gamma', 0\}$. Then the operator
    \begin{equation*}
        P \ede \rho_{0}^{\tilde \gamma}
        \rho_{\infty}^{\tilde \gamma'} (-\Delta + V)
        \in \Diff(\maS)\,, \quad \maS \ede \maSS\,.
    \end{equation*}
\end{proposition}

We let $\maS := \maSS$ from now on, with the parameters
$\tilde \gamma$ and $\tilde \gamma'$ as in the statement
of Proposition \ref{prop.inPsiS}.
We have $\Diff(\maS_{a, b, a', b'})\subset \Psi(\maS_{a, b, a', b'})$ (this is
valid for any CF-groupoid) and,
because the action of $\Psi(\maS_{a, b, a', b'})$ on
$\CIc(S^{n-1} \times [0, \infty])$ is injective, we see that
the operator $P := \rho_{0}^{\tilde \gamma}
\rho_{\infty}^{\tilde \gamma'} (-\Delta + V)$ corresponds to
a unique operator in $\Psi^{2}(\maS_{\tilde \gamma, \tilde \gamma-1,
\tilde \gamma' + 2, \tilde \gamma' + 1})$, which is elliptic 

by Equations (\ref{eq.schr3}-\ref{eq.schr6}). Hence it has a parametrix
$Q_{0} \in \Psi^{-2}(\maS_{\tilde \gamma, \tilde \gamma-1,
\tilde \gamma' + 2, \tilde \gamma' + 1})$. We obtain the
following result on the structure of the parametrices of
$-\Delta + V$.

\begin{theorem} \label{thm.inPsiS}
    There exists $Q \in \Psi^{-2}(\maS_{\tilde \gamma, \tilde \gamma-1,
    \tilde \gamma' + 2, \tilde \gamma' + 1})$ such that
    \begin{equation*}
        (- \Delta + V)\left(Q \rho_{0}^{\tilde \gamma}
        \rho_{\infty}^{\tilde \gamma'}\right)    \ , \ \
        \left(Q \rho_{0}^{\tilde \gamma}
        \rho_{\infty}^{\tilde \gamma'}\right) (- \Delta + V)
        \in \Psi^{-\infty} (\maS_{\tilde \gamma, \tilde \gamma-1,
    \tilde \gamma' + 2, \tilde \gamma' + 1}).
    \end{equation*}
\end{theorem}

Let $\phi := \rho_{0}^{\tilde \gamma} \rho_{\infty}^{\tilde \gamma'}$
(it is the same $\phi$ as the one appearing in the Introduction.)

\begin{proof}
    Let $Q := Q_{0}$ be the differential operators (parametrix)
    considered in the paragraph before the statement of
    the Theorem.

    The definition of $Q_{0}$ is
    such that
    \begin{equation*}
        \phi (- \Delta + V) Q
        \ , \ \
        Q \phi (- \Delta + V)
        \in \Psi^{-\infty}(\maS_{\tilde \gamma, \tilde \gamma-1,
        \tilde \gamma' + 2, \tilde \gamma' + 1})\,.
    \end{equation*}
    Then
    \begin{equation*}
        (- \Delta + V) \left( Q \phi \right) \seq
        \phi^{-1} \left[ \phi
        (- \Delta + V) Q \right] \phi
        \in \Psi^{-\infty}(\maS_{\tilde \gamma, \tilde \gamma-1,
        \tilde \gamma' + 2, \tilde \gamma' + 1})\,,
    \end{equation*}
    by Proposition \ref{prop.conj.rho}.
\end{proof}

Let us now sketch how one can treat the resolvents in our setting. Let 
$\maS := \maSS$, $\phi := \rho_{0}^{\tilde \gamma} \rho_{\infty}^{\tilde \gamma'}$, 
as before,  and let  $\overline{\Psi}^{-\infty}(\maS)$,
be the closure of $\Psi^{-\infty}(\maS)$ in the norms 
\begin{equation}\label{eq.norms}
    T \mapsto \|(\phi \Delta)^{i}T (\phi \Delta)^{j}\|_{\maL(L^{2}(\RR^{n}))}\,,
    \qquad i, j \in \ZZ_{+}\,,
\end{equation}

where $\maL(L^{2}(\RR^{n}))$ is the space of bounded operators on $L^{2}(\RR^{n})$
with the usual norm. Let $\overline{\Psi}^{m}(\maS) := \Psi^{m}(\maS) +
\overline{\Psi}^{-\infty}(\maS)$. Then

\begin{theorem}
    Let $z \in \CC$ and let us assume that $\phi(-\Delta + V) - z$ is invertible on
    $L^{2}(\RR^{n})$. Then $[\phi(-\Delta + V) - z]^{-1} \in
    \overline{\Psi}^{-2}(\maS).$
\end{theorem}

\begin{proof} 
    This follows from Proposition \ref{prop.inPsiS}
    and the results of \cite{LMNspectral, Vassout} (see 
    also \cite{SchroheFrechet}). Indeed,
    $\phi(-\Delta + V) - z \in \Diff(\maS) \subset \Psi(\maS) \subset
    \overline{\Psi}(\maS)$. The latter is stable under inverses, by the
    results of \cite{LMNspectral, Vassout}, hence $[\phi(-\Delta + V) - z]^{-1}
    \in \overline{\Psi}^{-2}(\maS)$, as claimed.
    
\end{proof}

Note that, while our statements about the structure of $-\Delta +V$
require only Lie-Rinehart algebras, the statement about its
parametrices does require a CF-groupoid integrating the relevant
Lie-Rinehart algebras.

\def\cprime{$'$}

\end{document}